\RequirePackage{fix-cm}

\documentclass[smallcondensed]{svjour3}     

\smartqed  
\usepackage{amsfonts}
\usepackage{amssymb}
\usepackage{enumitem}
\usepackage{mathrsfs}
\usepackage{tikz}
\usetikzlibrary{decorations.pathreplacing}
\usetikzlibrary{shapes.misc, positioning}
\usepackage{subfigure}
\usepackage{comment}
\usepackage{latexsym}
\usepackage{graphicx}
\usepackage{amsmath}
\usepackage{float}
\usepackage{subfig}
\usepackage{algorithm}
\usepackage{algorithmic}
\usepackage{multirow}
\usepackage{mathrsfs}
\usepackage{xcolor}
\usepackage[misc]{ifsym}

\makeatletter
\def\hlinew#1{%
	\noalign{\ifnum0=`}\fi\hrule \@height #1 \futurelet
	\reserved@a\@xhline}
\makeatother

\DeclareMathOperator{\diag}{diag}
\DeclareMathOperator{\Trace}{Trace}

\newcommand*{\affaddr}[1]{#1}
\newcommand*{\affmark}[1][*]{\textsuperscript{#1}}

\begin{document}

\title{Two-Grid based Adaptive Proper Orthogonal Decomposition Method for Time Dependent Partial Differential Equations
\thanks{This work was supported by  the National Key Research and Development Program of China under grant 2019YFA0709601,  the National Natural  Science Foundation of China under grants 91730302
and 11671389, the Key Research Program of Frontier Sciences of the Chinese Academy of Sciences
under grant QYZDJ-SSW-SYS010, and the NSF grant IIS-1632935.}
}
 \author{Xiaoying Dai\affmark[1] \and Xiong Kuang\affmark[1] \and \\ Jack Xin\affmark[2]\and Aihui Zhou\affmark[1]}
 \authorrunning{
	Xiaoying Dai \and
	Xiong Kuang \and
	Jack Xin \and
	Aihui Zhou
}
\institute{\Letter \quad  Xiaoying Dai \\	
	\hspace*{0.5cm} daixy@lsec.cc.ac.cn  \\
	 \affaddr{\affmark[1] \quad LSEC, Institute of Computational Mathematics and Scientific/Engineering Computing,
	 	Academy of Mathematics and Systems Science, Chinese Academy of Sciences,  Beijing 100190, China; and School of Mathematical Sciences,
	 	University of Chinese Academy of Sciences, Beijing 100049, China}\\
	 \affaddr{\affmark[2]    \quad Department of Mathematics, University of California at Irvine, Irvine, CA 92697, USA}\\
}
\date{Received: date / Accepted: date}
\maketitle
\begin{abstract}
	
In this article, we propose a two-grid based adaptive proper orthogonal decomposition (POD) method to solve the time dependent partial differential equations. Based on the error obtained in the coarse grid, we propose an error indicator for the numerical solution obtained in the fine grid. Our new  algorithm is cheap and  easy to be implement.  We apply  our new method to the solution of   time-dependent advection-diffusion equations with the Kolmogorov flow and the ABC flow. The numerical results show that our method is more efficient than the existing POD methods.
 \end{abstract}

\section{Introduction}
\label{intro}
Time dependent partial differential equations play an important role in scientific and engineering computing. Many physical phenomena  are described by time dependent partial differential equations, for example, the seawater intrusion \cite{bakker2000simple}, the heat transfer \cite{cannon1984one}, fluid equations \cite{biferale1995eddy,Xin-KPP}.  The design and analysis of high efficiency numerical schemes  for time dependent  partial differential equations has always been an active research topic.

For the spatial discretization of the time dependent partial differential equations, some classical discretization methods, for example, the finite element method \cite{brenner2007mathematical}, the finite difference method \cite{leveque2007finite}, the plane wave method \cite{kosloff1983fourier}, can be used. However, for many complex systems, these classical  discretization methods will usually result in  discretized systems with millons or even billions of degree of freedom, especially  when the spatial dimension is equal to or larger than three. 
 Therefore, if we  use these classical  discretization methods to deal with the spatial discretization at each time interval, the computational cost will be very huge \cite{bieterman1982finite,bueno2014fourier,smith1985numerical}.
 
We realize that many model reduction methods have been developed to reduce the degrees of freedom and then the  computational costs \cite{Reduced-Order-Methods}. The general idea for these methods is to project the orginal system onto a low-dimensional approximation subspace  so that the number of freedoms involved in the discretization will be  reduced significantly.
The proper orthogonal decomposition (POD) \cite{lumley-1967,sirovich-1987} method is one of the most commonly used ways to define the low-dimensional subspace. Other 
model reduction methods include reduced basis methods\cite{boyaval2010reduced,hesthaven-rozza-2015,maday2006reduced,maday2002reduced}, balanced truncation 
\cite{gugercin2004survey} and Krylov subspace methods \cite{feldmann1995efficient}. We refer to  \cite{benner-gugercin-willcox-2015,chinesta-huerta-willcox-2017} for the reviews of the model reduction methods.

 The POD method has been widely used in  time dependent partial differential equations, see, e.g.,  \cite{atwell2001proper,benner-gugercin-willcox-2015,Burgers-equation,POD-arbitrary-FEM,POD-parabolic,APOD-Spectral-Xin,pinnau2008model}.  Some review about  POD method and its applications can be found in \cite{hesthaven-rozza-2015,volkwein-2011}. 
The basic idea of the POD method
 is to start with an ensemble of data, called snapshots, 
then construct POD modes by performing singular value decomposition 
(SVD) to these snapshots 
\cite{burkardt2006pod,ito1998reduced,ly2001modeling,pinnau2008model,sirovich-1987}. 
 After that, by projecting the orginal system onto the reduced-order subspace spanned by these POD modes, an POD reduced model is then obtained and solved instead of the orignal system.  By choosing the snapshots properly,  the number of the POD modes will usually be much less than the degrees of freedom resulted by the traditional methods. Therefore,  it will be much cheaper to discretize  the time dependent partial differential equations in the subspace spanned by these POD modes. 

We see that the choice of the snapshot ensemble is a crucial issue in constructing  POD modes, and  can greatly affect the approximation of solution for the original systems.  We refer to  \cite{kunisch-volkwein-2010} for a discussion on  this issue, which will not be addressed here.  

For the classical POD method for time dependent partial differential equations,  the snapshots are usually collected from trajectories of the dynamical system at certain discrete time instances  obtained by one of the traditional methods over some interval $[0,T_{0}]$ \cite{burkardt2006pod,sirovich-1987}.
 Once the POD modes are chosen, they will be fixed and not updated during the time evolution. However,  as the time evolution of the solution takes place, the feature of the solution will usually change. Therefore, if the POD modes are not updated any more, the approximation error obtained by the POD methods may become nonnegligible.

To improve the efficiency of the  POD methods, some adaptive POD approaches have been introduced in the past years. In \cite{peherstorfer-willcox15a}, the authors proposed dynamic reduced models for dynamic data-driven application
systems. The dynamic reduced models take into account the changes
in the underlying system by exploiting the opportunity presented by the dynamic sensor data and adaptively incorporating sensor data during the dynamic process. In \cite{APOD-two-Gs}, some local POD and Galerkin projection method was proposed, however,  the assumption that the Galerkin
approximation converges to the right solution if an appropriately large number of modes is used, which may not be true, especially for the case that the POD modes are obtained by the snapshots in some time subinterval other than the whole time interval.  In \cite{zhuang,APOD-residual-estimates,terragni2015simulation}, some adaptive POD methods for time dependent 
paritial equations are proposed by constructing some residual type error indicators to determine whether the POD modes need to be updated. The residual type  error indicator is efficient, however, since the calculation of the residual involves the data in the original high dimensional space,  it is usually expensive to compute.

In this paper, 
we propose a two-grid based adaptive POD method, 
where we use the error obtained in the coarse mesh to construct the error indicator to tell us if we need to update the POD subspace in the fine mesh or not. Since the degree of freedom corresponding to  the coarse mesh is much less than that corresponding to the fine mesh, it is cheap to calculate the error indicator. Therefore, by our method, we can easily compute the error indicator, and then update the POD subspace when needed.

We note that the two-grid approach has been  widely used in the discretization of different kinds of  partial differential equations \cite{dai-zhou-2008,xu-zhou-2000,xu-zhou-2001}. The main difference lies in that the coarse grid in this paper is only used to obtain the indicator which tells us at which instants we need to update the POD modes.   From some special perspective, our two-grid approach here can be viewed as one of the widely used multi-fidelity type approaches. The multi-fidelity methods accelerate the
solution of original systems by combining high-fidelity and low-fidelity model evaluations, and have been widely used in many different applications \cite{forrester-sobester-keane-2007,ng-willcox-2014,teckentrup-jantsch-webster-2015}.  Some review about multi-fidelity methods can be found in \cite{peherstorfer-willcox-gunzburger-2018}.

The rest of this paper is organized as follows. First, we give some preliminaries in Section 2, including a brief  introduction for the finite element method, the standard POD-Galerkin method, the general framework for adaptive POD method,  and the residual based adaptive POD method. Then, we propose our two-grid based adaptive POD method in Section 3. Next, we apply our new method to the simulation of some typical time dependent partial differential equations, including advection-diffusion equation with three-space dimensional velocity field, such as Kolmogorov flow and ABC flow, and use these numerical experiments to show the efficiency  and the advantage of our method to the existing methods in Section 4. Finally, we give some concluding remarks in Section 5, and provide some additional numerical  experiments for tuning the parameters and the coarse mesh size in Appendix A and Appendix B, respectively. 

\section{Preliminaries}
\setcounter{equation}{0}
We consider the following general time dependent partial differential equation
\begin{eqnarray}\label{question1}
\left\{
\begin{aligned}
&u_{t}- D_{0} \Delta u +\vec B(x,y,z,t) \cdot \nabla u + c(x,y,z,t)u = f(x,y,z,t) ,\quad  \text{in} \Omega \times (0, T]\\
&u(x,y,z,0)=h(x,y,z),\\
&u(x+l,y,z,t)=u(x,y+l,z,t)=u(x,y,z+l,t)=u(x,y,z,t),\\
\end{aligned}
\right.
\end{eqnarray}
where $\Omega = [0, l]^3, f \in L^2(0,T; L^2(\Omega)), c\in L^{\infty}(\Omega), \vec{B} \in C(0,T;W^{1,\infty}(\Omega)^3)$ and $D_0$ is a constant. \\
\indent Define a bilinear form
\begin{eqnarray*}\label{bilinear}
	a(t; u, v) = D_0 (\nabla u,\nabla v)-D_0\int_{\partial \Omega} \frac{\partial u}{\partial n}vd\sigma + (\vec{B} \cdot \nabla u, v)+(cu,v), \forall u, v \in H^1(\Omega),
\end{eqnarray*}
where $(\cdot , \cdot)$ stands for the inner product in $L^2(\Omega)$, and the function space
\begin{eqnarray*}
	V = \{v \in H^1(\Omega):v|_{x=0} = v|_{x=l},v|_{y=0} = v|_{y=l}, v|_{z=0} = v|_{z=l}\}.
\end{eqnarray*}
Then the variational form of equation (\ref{question1}) can be written as follows: find $u \in V$ such that
\begin{eqnarray}\label{original-variational-form}
(\frac{\partial u }{\partial t}, v)+a(t; u, v)=  (f(x, y, z, t), v),\forall v\in V.
\end{eqnarray}

In order to solve (\ref{original-variational-form}) numerically, we choose the implicit Euler scheme  \cite{acary2010implicit,tone2006long} for the temporal discretization. We partition the time interval into $N \in \mathbb{N}$ subintervals with equal length $\delta t = T\slash N$, and set $u^k(x,y,z) = u(x, y, z, t_k)$ where $t_k = k*\delta t$, for $k \in \{0,1,\ldots,N\}$.
Then the semi-discretization scheme of (\ref{original-variational-form}) can be written as:
\begin{eqnarray}\label{finite-element-time}
\left.
\begin{aligned}
(\frac{u^k(x,y,z)-u^{k-1}(x,y,z)}{\delta t},v)+ a(t_k; u^k(x,y,z), v) = (f(x, y, z, t_k), v),\forall v\in V.\quad
\end{aligned}
\right.
\end{eqnarray}
 
\subsection{Standard finite element method}\label{sec:2}
In this subsection, we introduce the standard finite element discretization for the  equation (\ref{finite-element-time}) briefly. For more detailed introduction on standard finite element method, please refer to e.g. \cite{shen2013finite,thomee1984galerkin}. \\
\indent Let $\mathcal{T}_h$ be a regular mesh over $\Omega$,  that is, there exists a constant $\gamma^{\ast }$ such that \cite{chen2014adaptive}
\begin{eqnarray*}
	\frac{h_\tau}{\rho_\tau} < \gamma^{\ast}, \forall \tau \in \mathcal{T}_h
\end{eqnarray*}
where $h_\tau$ is the diameter of $\tau$ for each $\tau \in \mathcal{T}_h $ and  $\rho_{\tau}$ is the diameter of the biggest ball contained in $\tau \in \mathcal{T}_h,  h = \max{h_\tau, \tau \in \mathcal{T}_h}$. Denote $\# \mathcal{T}_h$ the  degree of freedom of mesh $\mathcal{T}_h$. Define the finite element space as
$$V_{h}=\{v_h:v_h|_e \in P_{e}, \forall e \in \mathcal{T}_h  \quad \text{and} \quad  v_h \in C^{0}(\bar{\Omega})\}\cap V,$$
where $P_e$ is a set of polynomial function on element $e$.

Let $\{ \phi_{h,i} \}_{i=1}^n$ be a basis for $V_h$, that is
$$V_h := \text{span} \{\phi_{h,1}, \phi_{h,2},\ldots, \phi_{h,n} \}.$$
Then the numerical approximation of $u^k(x,y,z)$ can be expressed as
\begin{equation}\label{FEM-form}
u^k_h(x,y,z) = \sum \limits_{i=1}^{n} \beta^k_{h,i}\phi_{h,i}(x,y,z).
\end{equation}
Inserting (\ref{FEM-form})
into (\ref{finite-element-time}), and setting $v= \phi_{h,j}, j=1,2,3,\ldots,n$, respectively, we obtain
\begin{eqnarray}\label{finite-element-space}
\left.
\begin{aligned}
\sum\limits_{i=1}^{n}[\frac{\beta^k_{h,i}-\beta^{k-1}_{h,i}}{\delta t}(\phi_{h,i},\phi_{h,j} )+ \beta^{k}_{h,i} a(t_k; \phi_{h,i}, \phi_{h,j})] = (f(x, y, z, t_k), \phi_{h,j}).
\end{aligned}
\right.
\end{eqnarray}
We can rewrite (\ref{finite-element-space}) as
\begin{eqnarray}\label{No-Matrix-form}
\left.
\begin{aligned}
\sum\limits_{i=1}^{n} \beta^k_{h,i} [(\phi_{h,i},\phi_{h,j} ) +\delta t  a(t_k; \phi_{h,i}, \phi_{h,j})] =\delta t (f(x, y, z, t_k), \phi_{h,j})\\ +\sum\limits_{i=1}^{n}\beta^{k-1}_{h,i} (\phi_{h,i}, \phi_{h,j}).\qquad \qquad \qquad
\end{aligned}
\right.
\end{eqnarray}
Define
 \begin{equation*}
\begin{aligned}
& \vec{A}^k_{h,ij} =(\phi_{h,j}, \phi_{h,i}) + \delta t a(t_k; \phi_{h,j}, \phi_{h,i}) ,\quad \vec{u}^{k}_h=(\beta^k_{h,1}, \beta^k_{h,2},\beta^k_{h,3},\ldots,\beta^k_{h,n})^{T}, \\
&	\vec{b}_h^k  = \delta t*((f, \phi_{h,1}),\ldots,(f, \phi_{h,n}))^{T}, \quad \vec{C}_{h,ij} = (\phi_{h,j}, \phi_{h,i}).
\end{aligned}
\end{equation*}
Then (\ref{No-Matrix-form}) can be written as the following algebraic form
\begin{equation}\label{FEM-equation}
\vec{A}^k_{h} \vec{u}^{k}_h = \vec{b}^k_h + \vec{C}_h \vec{u}^{k-1}_h.
\end{equation}
\subsection{POD method} \label{subsec-POD}

\begin{figure}
	\centering
	\begin{tikzpicture}
	\coordinate [label=above: $0$] (A1) at (0, 0);
	\coordinate [label=above:   \footnotesize $T_0$] (B1) at (3, 0);
	\draw  (A1) -- (B1);
	\draw [decorate,decoration={brace, mirror, amplitude=5pt},xshift=0 pt,yshift= 0pt]
	(A1) -- (B1) node [black, midway, yshift= -0.5cm] {\footnotesize $I_{\mbox{FEM}}$};
	
	\coordinate [label=above:\footnotesize $T$] (C1) at (11, 0);
	\draw  (B1) -- (C1);
	\draw [decorate,decoration={brace, mirror, amplitude=5pt},xshift=0 pt,yshift= 0pt]
	(B1) -- (C1) node [black, midway, yshift= -0.5cm] {\footnotesize $I_{\mbox{POD}}$};
 	\end{tikzpicture}
	\caption{Sketch of the POD method. $I_{\mbox{FEM}}$ and $I_{\mbox{POD}}$ refer to the time intervals where (\ref{finite-element-time}) is discretized  in $V_h$ and the subspace spanned by the POD modes, respectively.}\label{fig-POD}
\end{figure}
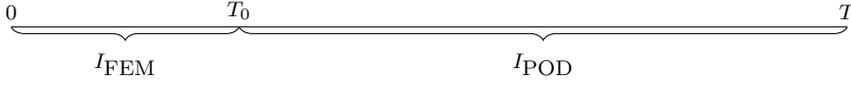

For the classical POD method for time dependent partial differential equations, the snapshots are usually collected from trajectories of the dynamical system at certain discrete time instances  obtained by one of the traditional methods over some interval $[0,T_{0}]$. Then, the POD modes are constructed  by performing SVD  to these snapshots.  Projecting the original dynamic system onto the reduced-order subspace spanned by the POD modes, we then obtain the POD reduced model. To understand this process clearly, we show the sketch of the POD method in Fig. \ref{fig-POD}.

 Now, we turn to introduce the details for each step of POD method for  discreting problem (\ref{finite-element-time}).  
\begin{itemize}[leftmargin=*]
	\item[1.] {\bf Snapshots} \\
	Discretize (\ref{finite-element-time}) in $V_h$ on the interval $[0,T_{0}]$, and collect the numerical solution per $\delta M$ steps($\delta M$ is a parameter to be specified  in the numerical experiment). Set $n_s =\lfloor \frac{T_0}{\delta t \cdot  \delta M} \rfloor $, where $\lfloor * \rfloor$ means the round down, and denote
	$$\vec{U}_h = [\vec{u}^0_h, \vec{u}^{\delta M}_h, \ldots,\vec{u}^{n_s \cdot \delta M}_h],$$
 	\item[2.] {\bf POD modes}\\ 
 	Perform SVD to the snapshots matrix $\vec{U}_h \in \mathbb{R}^{n\times (n_s+1)}$, and obtain
	\begin{equation}\label{SVD1}
			\vec{U}_h =\vec{R} \vec{S} \vec{V}^T,
	\end{equation}
 	where $\vec{R}= \left [R_1, R_2,\ldots, R_r \right ] \in \mathbb{R}^{n \times r}$, $\vec{V} = [V_1, \ldots, V_r] \in \mathbb{R}^{(n_s+1) \times r}$  are the left and right projection matrices,
	respectively, and $\vec{S} = \diag\{\sigma_1, \sigma_2,\cdots, \sigma_{r}\}$
	with $ \sigma_1 \geq \sigma_2 \geq\cdots \geq\sigma_{r} > 0$. The rank of $\vec{U}_h$ is  $r$, and obviously $ r \leq \min(n, n_s+1)$.
	
	\qquad Set $m = \min \{k| \sum\limits_{i=1}^{k}\vec{S}_{i, i}>\gamma_1*  \Trace(\vec{S})\}$($\gamma_1$ is a parameter to be specified  in the numerical experiment), then the POD modes  are constructed by
	\begin{equation}
	(\psi_{h,1}, \psi_{h,2},\ldots,\psi_{h,m})=(\phi_{h,1}, \phi_{h,2}, \ldots,\phi_{h,n}) \vec{\widetilde{R}},
	\end{equation}
	where $\vec{\widetilde{R}}=[R_1,R_2,\ldots,R_{m}]$.
	
   \qquad For the convenience of the following discussion, we summarize the process of constructing POD modes as routine POD$\underline{\hbox to 0.2cm{}}$Mode$(\vec{U}_h, \gamma_1, \Phi_h, m, \Psi_h)$, where $\Phi_h=(\phi_{h,1}, \phi_{h,2}, \ldots,\phi_{h,n})$ and $\Psi_h=(\psi_{h,1}, \psi_{h,2}, \ldots,\psi_{h,m})$, please see  Algorithm \ref{construc-POD} for the details.

	\begin{minipage}{11cm}
		\begin{algorithm}[H]
			\caption{  POD\underline{\hbox to 0.2cm{}}Mode$(\vec{U}_h, \gamma, \Phi_h, m, \Psi_h)$}\label{construc-POD}
			\begin{algorithmic}
				\REQUIRE
				$\vec{U}_h$,$\gamma$, $\Phi_h=(\phi_{h,1}, \phi_{h,2},  \ldots,\phi_{h,n})$.
				\ENSURE
				$m$ and  POD modes $\{\psi_{h,1}, \psi_{h,2},\ldots,\psi_{h,m}\}$.
				\STATE{Step1: Perform SVD to $\vec{U}_h$, and obtain $\vec{U}_h = \vec{R}\vec{S}\vec{V}^T$, where $\vec{S} = diag(\sigma_1, \sigma_2,\ldots, \sigma_r)$ with $\sigma_1\geq\sigma_2\geq \cdots \geq \sigma_r>0$.}
				\STATE{Step2: Set $m=\min \{k|\sum\limits_{i=1}^{k} \vec{S}_{i, i}>\gamma*  \Trace(\vec{S})\}$.}
				\STATE Step3: $(\psi_{h,1}, \psi_{h,2},\ldots,\psi_{h,m})=\Phi_h \vec{R}[:,1:m]$.\\
			\end{algorithmic}
		\end{algorithm}
	\end{minipage}
	\\\\
	\item[3.] {\bf Galerkin projection}\\
	\qquad  For $t>T_0$,  we discretize (\ref{finite-element-time}) in   $\mbox{span}\{\psi_{h,1}, \psi_{h,2},\ldots,\psi_{h,m}\}$, which is also called the POD subspace. That is, the solution of (\ref{finite-element-time}) is approximated by
	\begin{equation}\label{PODsolution-formule}
	u^k_{h,\mbox{POD}}(x,y,z) =\sum\limits_{i=1}^{m}\tilde{\beta}^k_{h,i} \psi_{h,i}(x,y,z).
	\end{equation}
	Inserting (\ref{PODsolution-formule})
	into (\ref{finite-element-time}), and setting $v= \psi_{h,j}, j=1, 2,3,\ldots,m $, respectively, we then obtain the following discretized problem:
 	\begin{equation}\label{PODtmp-equation}
	\begin{aligned}
	\sum\limits_{i=1}^{m} \tilde{\beta}^k_{h,i}[(\psi_{h,i},\psi_{h,j}) +\delta t  a(t_k; \psi_{h,i}, \psi_{h,j})] =\delta t (f(x, y, z, t_k), \psi_{h,j})\\ +\sum\limits_{i=1}^{m}\tilde{\beta}^{k-1}_{h,i}(\psi_{h,i}, \psi_{h,j}).
	\end{aligned}
	\end{equation}
	We define
\begin{equation*}
\begin{aligned}
& \vec{\widetilde{A}}^k_{h,ij} =(\psi_{h,j}, \psi_{h,i}) + \delta t a(t_k; \psi_{h,j}, \psi_{h,i}),  \\
&	\vec{\tilde{b}}_h^k  = \delta t*((f, \psi_{h,1}),\ldots,(f, \psi_{h,m}))^{T}, \quad \vec{\widetilde{C}}_{h,ij} = (\psi_{h,j}, \psi_{h,i}),\\
& \quad \vec{u}_{h,\mbox{POD}}^{k}=(\tilde{\beta}^k_{h,1}, \tilde{\beta}^k_{h,2},\tilde{\beta}^k_{h,3},\ldots,\tilde{\beta}^k_{h,m})^{T}.
\end{aligned}
\end{equation*}	Then (\ref{PODtmp-equation}) can be written as the following algebraic form
 	\begin{equation}\label{matrix_POD system}
	\vec{\widetilde{A}}^k_h\vec{u}_{h,\mbox{POD}}^{k} = \vec{\widetilde{b}}^k_h +\vec{\widetilde{C}}_h\vec{u}_{h,\mbox{POD}}^{k-1}. 
	\end{equation}
 \end{itemize}

By some simple calculation, we have that
	\begin{eqnarray*}
		\left.
		\begin{aligned}
			&\vec{\widetilde{A}}^k_h = \vec{\widetilde{R}}^{T}\vec{A}_h^k\vec{\widetilde{R}}, \quad
			\vec{\widetilde{b}}^k_h=  \vec{\widetilde{R}}^{T}*\vec{b}^k_h,\quad
			\vec{\widetilde{C}}_h = \vec{\widetilde{R}}^{T}\vec{C}_h\vec{\widetilde{R}}.\\
		\end{aligned}
		\right .
	\end{eqnarray*}

Summarizing the above discussion, we then get the standard POD  method  for discretizing (\ref{finite-element-time}), which is shown as Algorithm \ref{POD-Galerkin-method}.
\begin{algorithm}[H]
	\caption{POD method}\label{POD-Galerkin-method}
	\begin{algorithmic}[1]
		\STATE {Give $\delta t, \gamma_1, T_0$, $\delta M$, and the mesh $\mathcal{T}_h$. Set $n_s=\lfloor \frac{T_0}{\delta t \cdot  \delta M} \rfloor$. }
		\STATE {Discretize (\ref{finite-element-time}) in $V_h$ on interval $[0, T_0]$,  and obtain  $u^k_h,  k = 0, \cdots,  \lfloor T_0/\delta t \rfloor]$.}
		\STATE {Take snapshots $\vec{U}_h$ at   $t_0, t_{\delta M},\ldots, t_{n_s \delta M}$, respectively, that is\\ \quad \quad \quad \quad $\vec{U}_h = [\vec{u}^0_h, \vec{u}^{\delta M}_h, \ldots,\vec{u}^{ n_s \cdot \delta M}_h].$}
		\STATE {Construct POD modes $\Psi_h$ by   POD$\underline{\hbox to 0.2cm{}}$Mode$(\vec{U}_h, \gamma_1, \Phi_h, m, \Psi_h)$.}
        \STATE {$t = T_0$.}
		\WHILE {$t \le T$}
		\STATE $t=t+\delta t, k=k+1$.
		\STATE Discretize (\ref{finite-element-time}) in the subspace $\text{span} \{\psi_{h,1},\psi_{h,2},\ldots,\psi_{h,m}\}$, and  obtain $\vec{u}_{h,\mbox{POD}}^k$.
		\ENDWHILE
		\label{POD-algorithm}
	\end{algorithmic}
\end{algorithm}
\subsection{Adaptive POD method}\label{residual-APOD}
As the time evolution of the system, the feature of the solution may change a lot. Therefore, to improve the efficiency of the POD method, the adaptive updation of the POD mdoes are required.  Similar to the adaptive finite element method \cite{dai2008convergence}, the adaptive POD method consists of the loop of the form
\begin{center}
	\begin{equation}\label{loops}
	\begin{tabular}{|p{80mm}|}\hline
	\vspace{1mm}
	$
	~~~~~ \mbox{Solve}\rightarrow
	\mbox{Estimate}\rightarrow\mbox{Mark}\rightarrow\mbox{Update}.
	$ \vspace{2mm}\\
	\hline
	\end{tabular}
	\end{equation}
\end{center}
\vspace{2mm}

Here, we introduce this loop briefly. Suppose we have the POD modes \\
$\{\psi_{h, 1}, \cdots, \psi_{h, m}\}$, which can be written as
\begin{equation*}
	(\psi_{h,1}, \psi_{h,2},\ldots,\psi_{h,m})=(\phi_{h,1}, \phi_{h,2}, \ldots,\phi_{h,n}) \vec{\widetilde{R}},
\end{equation*}
	where $\vec{\widetilde{R}}=[R_1,R_2,\ldots,R_{m}]$.

First, we discretize the  equation (\ref{finite-element-time}) in the space spanned by these POD modes $\{\psi_{h, 1}, \cdots, \psi_{h, m}\}$,
 and obtain the POD approximation.
Then, we construct an error indicator to \textbf{Estimate} the error of the POD approximation. Next, we \textbf{Mark} the time instance where the  error indicator is too large. Finally, we  \textbf{Update} the POD modes at the marked time instance, and obtain the new POD modes. We repeat this loop in the new time interval until the terminal point.

The design of error indicator is an essential part in the \textbf{Estimate} step. For the error indicator, there are usually two requirements, one is that it can
estimate the error of the approximate solution very well, the other is that it should be very cheap to compute. In fact, the main difference between different adaptive POD methods lies in the construction of the error indicator. In this step, we design and calculate the error indicator $\eta_k$ at each time instance $t=k\delta t$.

In the \textbf{Mark} step, we mark the time instance $t_k$ if $\eta_k>\eta_0$ where $\eta_0$ is a user defined  threshold, which means  that the error of  POD approximation is so large that we need to update the POD modes.

For the \textbf{Update}  of the POD modes, we can do as follows. 
When we come to the time instance $t=p \delta t$ which  is marked, we go back to the previous time instance $t=(p-1)\delta t$. Starting from the time instance $t=(p-1)\delta t$,  we discretize  (\ref{finite-element-time}) in the finite element space $V_h$ and obtain $u^p_h, u^{p+1}_h,\ldots, u^{p+\lfloor \frac{\delta T}{\delta t} \rfloor}_h$. Set $n_{s_1} =\lfloor \frac{\delta T}{\delta t \cdot  \delta M} \rfloor $ and  denote
$$\vec{W}_{h,1} = [\vec{u}^p_h, \vec{u}^{p+\delta M}_h,\ldots,\vec{u}^{p+n_{s_1} \cdot \delta M}_h].$$
Performing SVD to $\vec{W}_{h,1}$,  we get $$\vec{W}_{h,1} = \vec{R}_1\vec{S}_1\vec{V}_1.$$
Here $\vec{R}_1 \in \mathbb{R}^{n \times r_1}$, $\vec{V}_1 \in \mathbb{R}^{(n_{s_1}+1) \times r_1}$  are the left and right projection matrices, respectively, and $\vec{S}_1 = \diag\{\sigma_{1,1}, \sigma_{1,2},\cdots, \sigma_{1, r_1}\}$ with $ \sigma_{1,1} \geq \sigma_{1, 2} \geq\cdots \geq\sigma_{1, r_1} > 0$. The rank of $\vec{W}_{h,1}$ is  $r_1 \leq \min(n, n_{s_1}+1)$.

Set $m_1=\min\{k:\sum\limits_{i=1}^k \vec{S}_{1, ii}\ge \gamma_2 \Trace(\vec{S}_1)\}$($\gamma_2$ is a parameter to be specified  in the numerical experiment),
then we combine the first $m_1$ column of $\vec{R}_1$ and  the old (previously used) POD modes, and get new matrix $\vec{W}_{h,2}$. That is
\begin{equation*}
\vec{W}_{h,2} = [\vec{R}_1[:,1:m_1], \vec{\tilde{R}}].   
\end{equation*}

We then preform   SVD to $ \vec{W}_{h,2}$ and obtain
\begin{equation}
\vec{W}_{h,2} = \vec{R}_2\vec{S}_2\vec{V}_2.
\end{equation}
Here $\vec{R}_2 \in \mathbb{R}^{n \times r_2}$ and $\vec{V}_2\in \mathbb{R}^{(m_1+m) \times r_2}$  are the left and right projection matrices, respectively, and $\vec{S}_2 = \diag\{\sigma_{2,1}, \sigma_{2,2},\cdots, \sigma_{2, r_2}\}$ with $ \sigma_{2,1} \geq \sigma_{2,2} \geq\cdots \geq\sigma_{2,r_2} > 0$. The rank of $\vec{W}_{h,2}$ is  $r_2$, $r_2 \leq \min(n, m_1+m)$.

Set $m_2=\min\{k:\sum\limits_{i=1}^k \vec{S}_{2, ii}\ge \gamma_3 \Trace(\vec{S}_2)\}$($\gamma_3$ is a parameter to be specified  in the numerical experiment), then the new POD modes are
\begin{equation}
\psi^{new}_{h,i} = \sum\limits_{j=1}^{n}\vec{R}_{2, ji}\phi_{h,j},\quad i=1,2,\ldots,m_2.
\end{equation}

For the convenience of the following discussion, we summarize the process of updating the POD modes as routine
Update\underline{\hbox to 0.2cm{}}POD\underline{\hbox to 0.2cm{}}Mode$(\vec{W}_{h,1},  \gamma_2, \gamma_3, \Phi_h, m, \Psi_h)$, which is shown as Algorithm \ref{updata-POD}.
\begin{algorithm}[H]
	\caption{Update$\underline{\hbox to 0.2cm{}}$POD$\underline{\hbox to 0.2cm{}}$Mode$(\vec{W}_{h,1},   \gamma_2, \gamma_3, \Phi_h, m, \Psi_h)$}\label{updata-POD}
	\begin{algorithmic}
		\REQUIRE
		$W_{h,1}$,$\gamma_2$,$\gamma_3$, $\Phi_h=(\phi_{h,1}, \phi_{h,2},  \ldots,\phi_{h,n})$, $m$ and $m$ old POD modes $\Psi_h$, $\Psi_h$ = $\Phi_h \vec{\tilde{R}}$.
		\ENSURE
		new $m$ and new $m$ POD modes $\{\psi_{h,1}, \psi_{h,2},\ldots,\psi_{h,m}\}$.
		\STATE{step1: preform SVD to $\vec{W}_{h,1}$, and obtain $\vec{W}_{h,1} = \vec{R}_1\vec{S}_1\vec{V}_1$, where $\vec{S}_1= \diag\{\sigma_{1,1}, \sigma_{1,2},\cdots, \sigma_{1, r_1}\}$ with $\sigma_{1,1} \geq \sigma_{1,2} \geq\cdots \geq\sigma_{1, r_1} > 0$. }
		\STATE{step2: set $m_1 = \min \{k| \sum\limits_{i=1}^{k}\vec{S}_{1,ii}>\gamma_2*  \Trace(\vec{S}_1)\}$.}
		\STATE{ step3: preform SVD to the matrix $\vec{W}_{h,2} = [\vec{R}_1[:,1:m_1], \vec{\tilde{R}}]$,    and obtain $\vec{W}_{h,2} = \vec{R}_2\vec{S}_2\vec{V}_2$, where $\vec{S}_2 =\diag\{\sigma_{2,1}, \sigma_{2,2},\cdots, \sigma_{2,r_2}\}$ with $\sigma_{2,1} \geq \sigma_{2,2} \geq\cdots \geq\sigma_{2, r_2} > 0$.}
		\STATE step4: set $m = \min \{k| \sum\limits_{i=1}^{k}\vec{S}_{2,ii}>\gamma_3*  \Trace(\vec{S}_2)\}.$
		\STATE step5:  $(\psi_{h,1}, \psi_{h,2},\ldots,\psi_{h,m})=\Phi_h \vec{R}_2[:,1:m].$ \qquad \qquad
	\end{algorithmic}
\end{algorithm}

 Then we obtain the following framework of the adaptive POD method for discretizing problem (\ref{finite-element-time}), see Algorithm \ref{APOD-algorithm} for the details.
\begin{algorithm}[!htbp]
	\caption{Adaptive POD method}\label{APOD-algorithm}
	\begin{algorithmic}[1]
		\STATE {Give $\delta t,T_0$, $\delta T, \gamma_1, \gamma_2, \gamma_3, \delta M$  and the mesh $\mathcal{T}_h$. Set $n_s = \lfloor \frac{T_0}{\delta t \cdot  \delta M} \rfloor$}
		\STATE {Discretize  (\ref{finite-element-time}) in $V_h$  on interval $[0,T_0]$ and obtain $u_h^k, \forall k \in [0, {\lfloor {T_0}/{\delta t}\rfloor}]$ , then take snapshots $\vec{U}_h$ at different times $t_0, t_{\delta M},\ldots, t_{n_s \cdot \delta M}$, that is\\ \quad \quad \quad \quad $\vec{U}_h = [\vec{u}^0_h, \vec{u}^{\delta M}_h, \ldots,\vec{u}^{n_s\cdot \delta M}_h]$.}
		\STATE {Construct POD modes $\Psi_h$ by POD$\underline{\hbox to 0.2cm{}}$Mode$(\vec{U}_h, \gamma_1, \Phi_h, m, \Psi_h)$.}
		\STATE {$t = T_0$.}
        \WHILE {$t \le T$}
		\STATE $t=t+\delta t$, $k=k+1$.
		\STATE Discretize (\ref{finite-element-time}) in the space $\text{span} \{\psi_{h,1},\ldots,\psi_{h,m}\}$, and obtain  $u_{h,\mbox{POD}}^k$, then compute error indicator $\eta_k$ by some strategy.
		\IF{ $\eta_k > \eta_0$}
		\STATE {$t=t-\delta t, k=k-1$.}
		\STATE {Discretize (\ref{finite-element-time}) in $V_h$ on interval $[t,t+\delta T]$ to get $u_h^{k+i}$, $i=1, \ldots, \frac{\delta T}{\delta t}$, from which to get snapshots $\vec{W}_{h,1}$, then update POD models $\Psi_h$ by  Update$\underline{\hbox to 0.2cm{}}$POD$\underline{\hbox to 0.2cm{}}$Mode$(\vec{W}_{h,1},  \gamma_2, \gamma_3, \Phi_h, m, \Psi_h)$}. $k=k+\frac{\delta T}{\delta t}$.
		\ENDIF
		\ENDWHILE
	\end{algorithmic}
\end{algorithm}
\subsection{Residual based adaptive POD method}
The main difference between different adaptive POD methods is the way to construct the error indicator.
Among the existing adaptive POD methods, the residual based adaptive POD method is the most widely used \cite{zhuang,APOD-residual-estimates}. The residual based adaptive POD algorithm uses the residual to construct the error indicator. Based on the residual corresponding to the POD approximation, people construct the error indicator $\eta_k$ at time instance $t=k\delta t$ as:
\begin{equation}\label{residual-APOD error indicator}
\eta_k = \frac{\|\vec{A}^k_h \vec{\widetilde{R}}\vec{u}_{h,\mbox{POD}}^{k}- \vec{b}^k_h  -   \vec{C}_h \vec{\widetilde{R}}\vec{u}_{h,\mbox{POD}}^{k-1}  \|_2}{ \|\vec{b}^k_h  +  \vec{C}_h \vec{\widetilde{R}}\vec{u}_{h,\mbox{POD}}^{k-1}\|_2}.
\end{equation}

\section{Two-grid based adaptive POD method}\label{TG-APOD}
For the residual based adaptive POD algorithm, we can see from (\ref{residual-APOD error indicator}) that in order to calculate the error indicator, we need to go back to the finite element space $V_h$ to  calculate the residual,  which is too expensive. Here, we propose a two-grid based adaptive POD algorithm (TG-APOD).
The main idea is to construct two finite element spaces, the coarse finite element space and the fine finite element space,
then construct the POD subspace in the fine finite element space, and use the coarse finite element space to construct the error indicator to tell us when the POD modes need to be updated.

We first construct a coarse partition $\mathcal{T}_H$ for the space domain $\Omega$  with  mesh size $H$ which is much larger than $h$  and a coarse partition for the time domain with time step $\Delta t$ which is much bigger than $\delta t$.  The finite element space corresponding to the partition $\mathcal{T}_H$ is denoted as  $V_H$.
In our following discussion,  the coarse mesh means the coarse spacial mesh size $H$ together with the coarse time step $\Delta t$, and the fine mesh means the fine spacial mesh size $h$ together with the fine time step $\delta t$. For simplicity, we require  that there exist some integers $M_1 \gg 1$ and $M_2 \gg 1$, such that $\Delta t = M_1 \delta t$ and $H = M_2 h$.

We first discretize the partial differential equation (\ref{original-variational-form}) by the same time discretized scheme as that used for obtaining (\ref{finite-element-time}),
that is, the implicit Euler scheme, with the coarse time step $\Delta t$, to obtained the following equation, which is similar as (\ref{finite-element-time}) but with a coarse time step $\Delta t$,
\begin{eqnarray}\label{eq-tmp}
\left.
\begin{aligned}
(\frac{u^k(x,y,z)-u^{k-1}(x,y,z)}{\Delta t},v)+ a(t_k; u^k(x,y,z), v) = (f(x, y, z, t_k), v),\forall v\in V.\quad
\end{aligned}
\right.
\end{eqnarray}
Then, we discretize the above equation in $V_H$ and obtain the finite element approximation $\vec{u}^k_H$. We then discretize (\ref{eq-tmp}) by the adaptive POD method
to obtain its adaptive POD approximation   $\vec{u}_{H, \mbox{POD}}^{k}$ with the error indicator $\eta_k$ being defined as
\begin{equation}\label{def-error-indicator}
 	\eta_k = \frac{\| \vec{u}^k_H -\vec{u}_{H, \mbox{POD}}^{k}\|_2}{\|\vec{u}^k_H\|_2}.
\end{equation}
 	
Now, we introduce each step of the loop (\ref{loops}) one by one for our TG-APOD algorithm.

\begin{itemize}[leftmargin=*]
	\item[1.] {\bf  Solve.}	Discretize the problem (\ref{finite-element-time}) in the subspace spanned by  the POD modes $\{\psi_{h,1}, \cdots, \psi_{h, m}\}$.

	\item [2.] {\bf  Estimate.} For $t=k \Delta t$, we  calculate the error indicator defined in (\ref{def-error-indicator}), from which we decide if we need to update the POD modes.
  	
	\item[3.] {\bf  Mark.}	Giving a threshold $\eta_0$, we set
\begin{equation}\label{stepsize}
\mbox{flag}_k =
\begin{cases}
 1, & \mbox{if
 $\eta_k > \eta_0$},\\
0, & \mbox{otherwise}.
\end{cases}
\end{equation}
  \item[4.] {\bf  Update.}
 If $\mbox{flag}_k = 1$,  then we move to $t_{k-1}$, and discretize the equation (\ref{finite-element-time})  in the fine finite element space $V_h$ on time interval $[t_{k-1}, t_{k-1} + \delta T]$, to get $u_h^{k+i}$, $i=1, \ldots, \frac{\delta T}{\delta t}$, from which to get snapshots $\vec{W}_{h,1}$, then update POD models $\Psi_h$ by  Update$\underline{\hbox to 0.2cm{}}$POD$\underline{\hbox to 0.2cm{}}$Mode$(\vec{W}_{h,1},  \gamma_2, \gamma_3, \Phi_h,m, \Psi_h)$.
\end{itemize}

 From the discussion above, we  see that the steps {\bf Estimate} and {\bf Mark}  only depend  on the coarse mesh. Therefore,  we can  get the set of all the marked time instances from the calculation on the coarse mesh.

 We summarize the discussion above and obtain our two-grid based adaptive POD method, and state it as Algorithm \ref{TG-APOD}.
\begin{algorithm}[!htbp]
	\caption{Two-grid based adaptive POD method}\label{TG-APOD}
	\begin{algorithmic}[1]
		\STATE {Give coarse mesh $\mathcal{T}_H$  with coarse time step $\Delta t$ and fine mesh $\mathcal{T}_h$ with fine time step $\delta t$,  $T_0$,  $\mathbb{S}=\{\}$,$\delta T$, $\gamma_1$,  $\gamma_2$, $\gamma_3$,  $\eta_0$, $\delta M$.}
		\STATE {Discretize (\ref{finite-element-time}) in $V_H$ on interval $[0,T_0]$, and take the snapshots $\vec{U}_H$.}
		\STATE {Construct POD modes  $\Psi_{H}$ by  POD$\underline{\hbox to 0.2cm{}}$Mode$(\vec{U}_H, \gamma_1, \Phi_{H}, m, \Psi_H)$. }
         \STATE {$t = T_0$.}
		\WHILE {$t \le T$}
		\STATE $t=t+\Delta t$.
		\STATE Discretize (\ref{finite-element-time}) in $V_H$ to get $u_H^k$, and discretize (\ref{finite-element-time}) in the space $\text{span}\{ \psi_{H,1},\ldots,\psi_{H, m}\}$ to get  $u_{H,\mbox{POD}}^k$, $k=\frac{t}{\Delta t}$,  then compute the error $\eta_k$ by (\ref{def-error-indicator}).
		\IF{ $\eta(t) > \eta_0$}
		\STATE $t = t-\Delta t.$
		\STATE $\mathbb{S}=\mathbb{S}\cup \{t\}.$
		\STATE  Discretize (\ref{finite-element-time}) in $V_H$  on interval $[t,t+\delta T]$, and take snapshots $\vec{W}_{H,1}$, then update POD modes $\Psi_{H}$ by  Update$\underline{\hbox to 0.2cm{}}$POD$\underline{\hbox to 0.2cm{}}$Mode$(\vec{W}_{H,1},   \gamma_2, \gamma_3, \Phi_{H}, m, \Psi_H).$
		\ENDIF
		\ENDWHILE
		\STATE {Discretize (\ref{finite-element-time}) in $V_h$ on interval $[0,T_0]$, and get the snapshots $\vec{U}_h$.}
		\STATE {Construct POD modes $\Psi_{h}$ by  POD$\underline{\hbox to 0.2cm{}}$Mode$(\vec{U}_h, \gamma_1, \Phi_{h}, m, \Psi_h)$.}
		\STATE {$t=T_0$.}
       \WHILE {$t \le T$.}
		\STATE $t=t+\delta t$.
		\STATE Discretize (\ref{finite-element-time}) in the space $\text{span}\{\psi_{h,1},\ldots,\psi_{h,m}\}$, and obtain $u_{h,\mbox{POD}}^k$, $k=\frac{t}{\delta t}$.
		\IF{ $t \in \mathbb{S}$}
		\STATE  Discretize (\ref{finite-element-time}) in $V_h$ on interval $[t,t+\delta T]$, and take snapshots $\vec{W}_{h,1}$, then update POD mode $\Psi_{h}$  by  Update$\underline{\hbox to 0.2cm{}}$POD$\underline{\hbox to 0.2cm{}}$Mode$(\vec{W}_{h,1}, \gamma_2, \gamma_3, \Phi_{h}, m, \Psi_h).$
		\ENDIF
		\ENDWHILE
	\end{algorithmic}
\end{algorithm}
\begin{remark}
For our two-grid based adaptive POD method, the steps {\bf Estimate} and {\bf Mark} are all carried out in the coarse mesh.  Since $\Delta t \gg \delta t$, the number of calculation for the error indicator is much smaller than that carried in the fine time interval. Besides,  $H \gg h$ means $\# \mathcal{T}_H \ll \# \mathcal{T}_h$, which implies that the cost for calculating  the error indicator $\eta_k$ is cheap.  These two facts make our two-grid based adaptive POD method much cheaper than the existing adaptive POD methods.
\end{remark}
\section{Numerical examples} \label{sec-numerical}
In this section, we apply our new method to two types of fluid  equations, the Kolmogorov flow and the ABC flow, which will show the efficiency of our two-grid based adaptive POD algorithm. For these two types of equations, we compare our new algorithm with the POD algorithm and the residual based adaptive POD algorithm. In our experiments, we use the standard finite element approximation corresponding to the fine mesh as the reference solution, and the relative  error of approximation obtained by the POD algorithm, or the residual based adaptive POD algorithm, or our new two-grid based adaptive POD algorithm is calculated as
\begin{equation}\label{def-exact-error}
\text{Error} = \frac{\| \vec{u}_h^k - \vec{u}_{h, *}^k\|_2}{\|\vec{u}_h^k\|_2},
\end{equation}
where $\vec{u}_h^k$ and $\vec{u}_{h, *}^{k}$ represent finite element approximations and different kinds of POD approximations at  $t_k$.

Our numerical experiments are carried out on the high performance computers of
the State Key Laboratory of Scientific and Engineering Computing, Chinese Academy of Sciences, and our code is written based on the  adaptive finite element software platform PHG \cite{phg}.

In our following statements for numerical experiments, we will sometimes use ``Residual-APOD" to denote the residual based adaptive POD method, and use ``TG-APOD" to denote our two-grid based adaptive POD method.

\subsection{Kolmogorov flow}
We first consider the following advection-dominated  equation with the cosines term (Kolmogorov flow) \cite{Kolmogorov-flow,galloway1992numerical,obukhov1983kolmogorov}, 
\begin{eqnarray}\label{Kolmogorov-equation}
\left\{
\begin{aligned}
&u_{t}-\epsilon \Delta u +\vec B(x,y,z,t) \cdot \nabla u = f(x,y,z,t), \quad (x, y, z)\in\Omega, t\in [0, T],\\
&u(x,y,z,0)= 0, \\
&u(x+2\pi,y,z,t)=u(x,y+2\pi,z,t)=u(x,y,z+2\pi,t)=u(x,y,z,t),  \\
\end{aligned}
\right.
\end{eqnarray}
where
\begin{eqnarray*}
	\vec B(x,y,z,t)
	&=& (\cos(y), \cos(z),\cos(x))+(\sin(z), \sin(x), \sin(y))\cos(t), \\
	f(x,y,z,t)&=&-\cos(y)-\sin(z)*\cos(t),\\
	\Omega = [0, 2\pi&]^3 &, T=100.
\end{eqnarray*}

For this example, we have tested $4$ different cases with  $\epsilon=0.5, 0.1$, $0.05$, and $0.01$, respectively.
We divide the $[0, 2\pi]^3$ into tetrahedrons to get the initial mesh containing $6$ elements, then refine the initial mesh $22$ times uniformly using  bisection to get our fine meshes.
We set $\delta t =0.005$,  and the type of the finite element basis is the piecewise linear function. For the three cases with  $\epsilon=0.5, 0.1$, and $0.05$, the parameters are chosen as $T_0=1.5, \delta T=1, \delta M=5$; for case   of $\epsilon=0.01$, the parameters are chosen as $T_0= 5, \delta T=3,\delta M=20$.
In the two-grid based adaptive POD algorithm, we refine the initial mesh $16$ times uniformly using bisection to obtain the coarse mesh, and the time step of the coarse mesh is set to $0.09$.   In fact, we have done some numerical experiments for testing   different coarse mesh size, all the numerical results can be found in Appendix B, from which we can see that this choice for the coarse mesh size is the best one taking into account both the accuracy and the cpu time.  We use $36$ processors for the simulation.

 Next, we will state how we choose the parameters $\gamma_i$ ($i=1,2,3$) and $\eta_0$. 
 \begin{itemize}
 	\item[1)] For the choice of the parameters $\gamma_i$ ($i=1,2,3$), we have done several numerical experiments in Appendix A for tuning the these parameters, from which  we  find that  the case of $\gamma_1=\gamma_2 = 0.999$, $\gamma_3 = 1.0-1.0\times 10^{-8}$ is a good choice. Therefore, in the following numerical experiments, we fix them as
 	$\gamma_1=\gamma_2 = 0.999$, $\gamma_3 = 1.0-1.0\times 10^{-8}$. 
 	
 	In fact, the above choice is reasonable.  Since the parameters $\gamma_1$ and $\gamma_2$ are both for extracting POD modes from the snapshots obtained from the standard finite element approximation of the original dynamic system, we do not need to set them too close to $1$. In fact, $\gamma_1$ and $\gamma_2$ being too close to $1$  may result in too many POD modes, which will cost too much cpu time. However, for $\gamma_3$, since it is used for extracting POD modes from the existing POD modes, we need to set to close to $1$, which means more information should be kept.
 	
 	\item[2)] For the choice of the parameter $\eta_0$, from the numerical experiments we have done and our  understanding for our two-grid based adaptive POD method,  we  observe that if the 
 	parameters  $\gamma_i$ ($i=1,2,3$) are  chosen properly, e.g., $\gamma_1=\gamma_2 = 0.999$, $\gamma_3 = 1.0-1.0\times 10^{-8}$, the magnitude of the error obtained by our method will be close to that of $\eta_0$, which implies that we do not need to tune the parameter $\eta_0$. In fact, if we need to obtain results with expected accuracy, we only need to set $\eta_0$ close to or a little smaller than that. In our following numerical experiments, we always choose $\eta_0 = 0.005$ for our two-grid based adaptive POD method.  
 \end{itemize}

 Anyway, in the following numerical experiments for our two-grid based adaptive POD method, we will choose  the parameters $\gamma_i$ ($i=1,2,3$) and $\eta_0$ as: 
 $\gamma_1=\gamma_2 = 0.999$, $\gamma_3 = 1.0-1.0\times 10^{-8}$, and $\eta_0 = 0.005$.

\begin{figure}[!htbp]	
	\begin{minipage}[t]{0.5\textwidth}
		\centering
		\includegraphics[width =\textwidth]{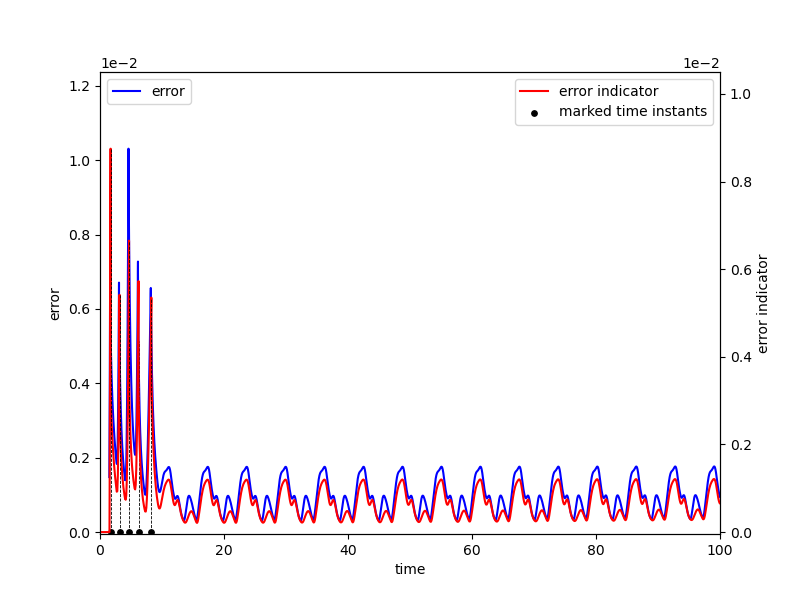}
		\caption*{a).$\epsilon=0.5$.}
	\end{minipage}
	\begin{minipage}[t]{0.5\textwidth}
		\centering
		\includegraphics[width = \textwidth]{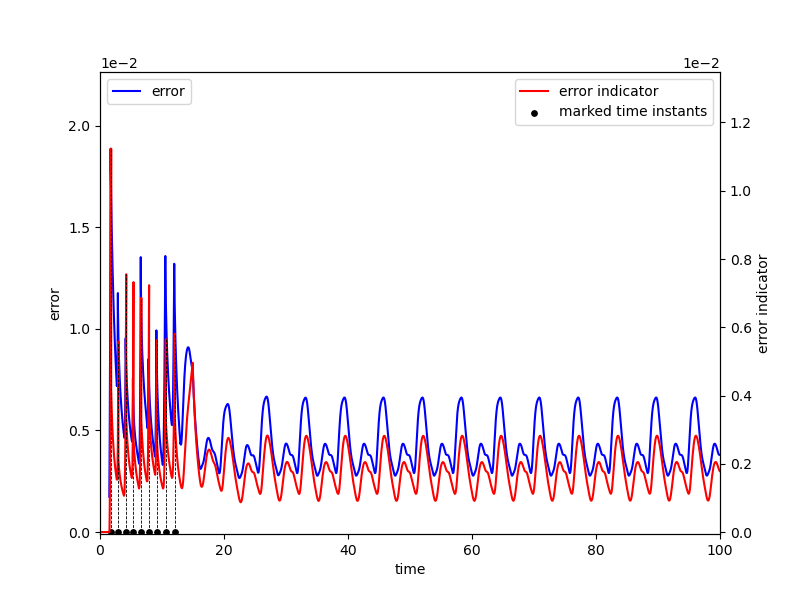}
		\caption*{b).$\epsilon=0.1$.}
	\end{minipage}
	\begin{minipage}[t]{0.5\textwidth}
		\centering
		\includegraphics[width = \textwidth]{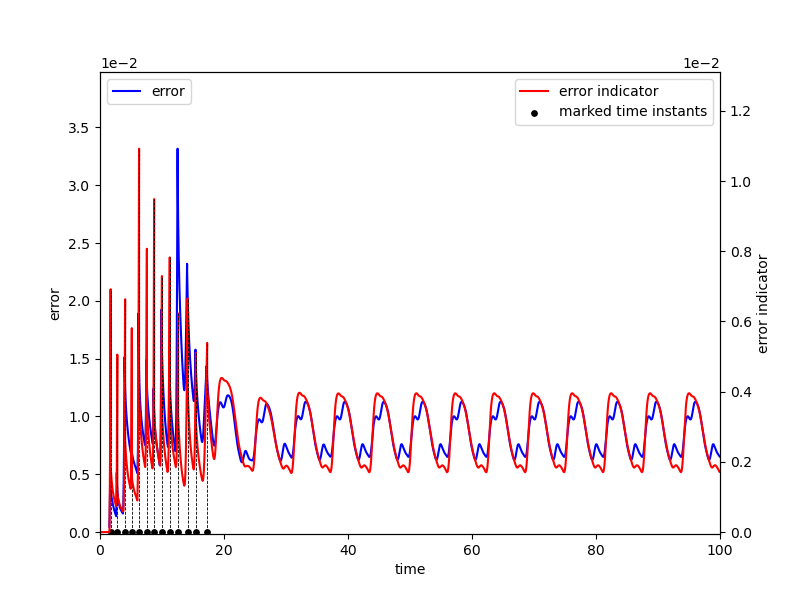}
		\caption*{c).$\epsilon=0.05$.}
	\end{minipage}	
	\begin{minipage}[t]{0.5\textwidth}
		\centering
		\includegraphics[width =\textwidth]{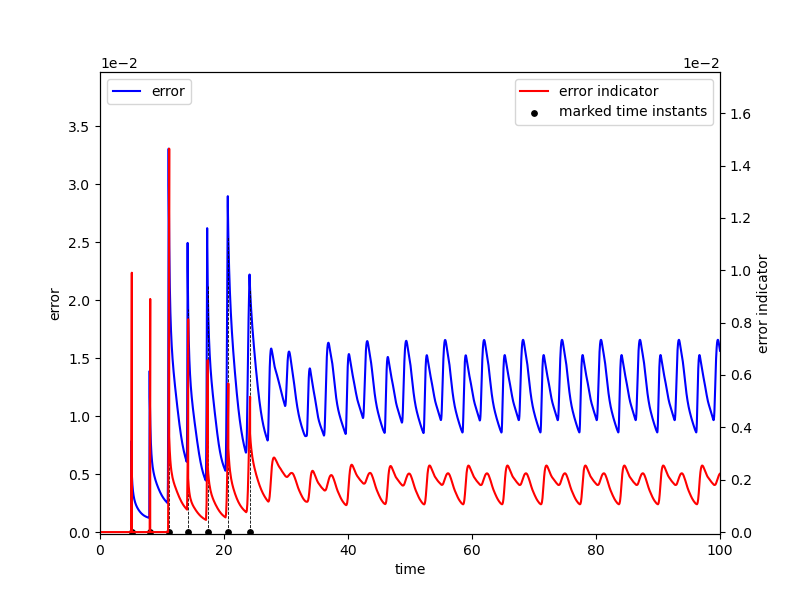}
		\caption*{d).$\epsilon=0.01$.}
	\end{minipage}
	\caption{The evolution curves of the error indicator and the error for solution of the Kolmogorov flow with $\epsilon=0.5, 0.1, 0.05,0.01$, respectively.}
	\label{err-errindicator1}
\end{figure}

We first use some numerical experiments to show the efficiency of our error indicator and the mark startegy. Fig. \ref{err-errindicator1} 
shows the evolution curves of  the error indicator and the error, and the time instance being marked to tell us when we should update the POD modes.

In Fig. \ref{err-errindicator1}, the y-axis in the left of each figure is the  error which is defined by (\ref{def-exact-error}), while the y-axis in the right of each figure is the error indicator defined by (\ref{def-error-indicator}), the black solid points denotes the time instances being marked which tell us when the POD modes need to be updated. 
	
By carefully comparing the curves of the error indicator and the error in Fig \ref{err-errindicator1}, we can see that curves of the errors and error indicators are similar, although their magnitudes are somewhat different. Note that in our two-grid based adaptive POD method, the error indicator is used to indicate the change trend of the error, which then tells us the time instance when the POD modes should be updated. By comparing the curves of the error and the time instances being marked to be updated, we observe that  the error indicator is efficient.

 We then compare the numerical results obtained by our TG-APOD method, the residual-APOD method, and the standard POD method. The detailed results are listed in Table \ref{compare3method-table1}.  The results obtained by standard finite element method are also provided as a reference. 

As we have stated before, for our  TG-APOD method, we choose   the parameters $\gamma_i$ ($i=1,2,3$) and $\eta_0$ as: $\gamma_1=\gamma_2 = 0.999$, $\gamma_3 = 1.0-1.0\times 10^{-8}$, and $\eta_0 = 0.005$. 
 While for the residual-APOD method, since people do not know how to set $\eta_0$ to obtained results with expected accuracy, they have to tune the parameter $\eta_0$. Therefore, to make the comparison reasonalble, we have tested the residual-APOD method by fixing   $\gamma_1=\gamma_2 = 0.999$, 
 $\gamma_3 = 1.0-1.0\times 10^{-8}$, but setting  different values for $\eta_0$ to find out the best performance of the residual-APOD method. 

In Table \ref{compare3method-table1}, `DOFs' means the degree of freedom, `Time' is the wall time for the simulation, `Average Error' is the average of error for numerical solution in each time instance. 

\begin{table}[!htbp]
	\centering
	\caption{The results of Kolmogorov flow with $\eta_0=0.005$ and  different $\epsilon$ obtained by Standard-FEM, POD, Residual-APOD, and TG-APOD, respectively.}
	\label{compare3method-table1}
	\begin{tabular}{|c|c|c|c|c|c|} \hlinew{1.5pt}
		$\epsilon$&Methods& $\eta_0$   &  DOFs& Average Error& Time(s)\\ \hlinew{1.2pt}
		& Standard-FEM & --& 4194304  & -- & 13786.49   \\
		0.5&	POD & -- &13 & 0.193343  & 547.75 \\
		&	Residual-APOD &$5\times 10^{-3}$ &9   & 0.151803  & 5590.62 \\
		&	Residual-APOD & $5\times 10^{-4}$ &  24 &  0.001166 &   6489.03 \\
		&	TG-APOD & $5\times 10^{-3}$ &24 &  0.001129  &1511.63\\ \hlinew{1.5pt} 	
		&	Standard-FEM &-- &4194304 & --& 11950.82  \\
		0.1&	POD &--&22  & 0.559466 & 535.47  \\
		&	Residual-APOD & $5\times 10^{-3}$ &28   & 0.167838 & 6241.83 \\
		&	Residual-APOD & $8\times 10^{-4} $   & 44  &  0.012444 & 7099.24 \\
		&	TG-APOD & $5\times 10^{-3}$& 54 & 0.004524 &  2693.41  \\ \hlinew{1.5pt}
		&	Standard-FEM &--&4194304& --  & 11379.64 \\
		0.05&	POD &-- &16 &  0.805803& 546.51 \\
		&	Residual-APOD & $5\times 10^{-3}$& 34  & 0.251224 &6470.52 \\
		&	Residual-APOD & $1\times 10^{-3}$   & 65  & 0.020072  &  7829.55\\
		&	TG-APOD & $5\times 10^{-3}$ &83 & 0.008663  & 3959.62  \\ \hlinew{1.5pt} 	
		&Standard-FEM &--&4194304 & --& 12955.44\\
		0.01&POD &--&43  & 0.771091 & 1265.95 \\
		&	Residual-APOD & $2\times 10^{-2}$  & 60  & 0.250940  & 7699.95 \\
		&Residual-APOD & $5\times 10^{-3}$& 97 & 0.051444 & 9182.46 \\
		&TG-APOD & $5\times 10^{-3}$&138  &   0.010969  & 6299.70 \\ \hlinew{1.5pt}
	\end{tabular}
\end{table}

We take a detailed look at Table  \ref{compare3method-table1}. First, we can easily find that the degrees of freedom for all the POD type methods are indeed much smaller than those for the standard finite element method, which means that these POD type methods have good performance in dimensional reduction. 
Then, we can see from  this example that the error of the numerical solution obtained by the standard POD method
 is too large, which makes the results meaningless. While both the residual based adaptive POD method and our two-grid based adaptive POD method can obtain numerical solution with higher accuracy,  which validates the effectiveness of the adaptive POD methods. Next, let's take a look at the results obtained by residual based adaptive POD method with diffferent $\eta_0$. We can see that as the decrease of the parameter $\eta_0$, the error of solutions obtained decreases while the cpu time cost increases.  For different $\epsilon$, the gap between the error obtained by the residual based adaptive POD method and $\eta_0$ are different. While for our two-grid based adaptive POD method, the magnitude of the error is close to that of $\eta_0$, even for $\epsilon = 0.01$. These coincides with what we have stated before. We also observe that as the decrease of $\epsilon$, the accuracy obtained by all the adaptive POD methods decreases, including our two-grid based POD method. We will pay more effort on the case with small $\epsilon$ in our future work.  
 For the comparison of our method with the residual based adaptive POD method, we  see that, by using the same parameters, the accuracy of results obtained by our method is much higher than those obtained  by the residual based adaptive POD method, while the cpu time cost by our method is much less. Although by choosing proper $\eta_0$, the residual based adaptive POD method can obtain results as accurate as our method,   the cpu time cost will be much more than that cost by our method. 
 Anyway, taking  both the accuracy and the cpu-time into account,  we can see clearly that our method is more efficient than the residual based adaptive POD method.

\begin{figure}[!htbp]
	\centering	
	\subfigure{
		\begin{minipage}[t]{0.45\linewidth}
			\centering
			\includegraphics[width=6cm,height=4cm]{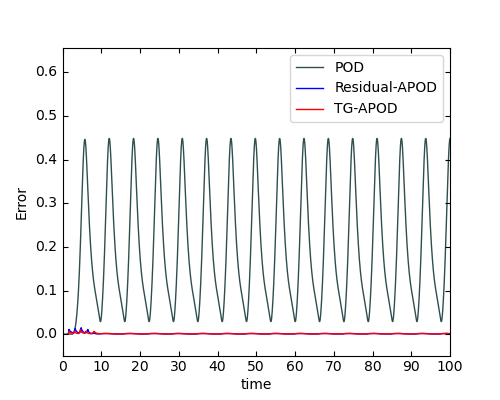}	
		\end{minipage}%
	}%
	\subfigure{
		\begin{minipage}[t]{0.45\linewidth}
			\centering
			\includegraphics[width=6cm,height=4cm]{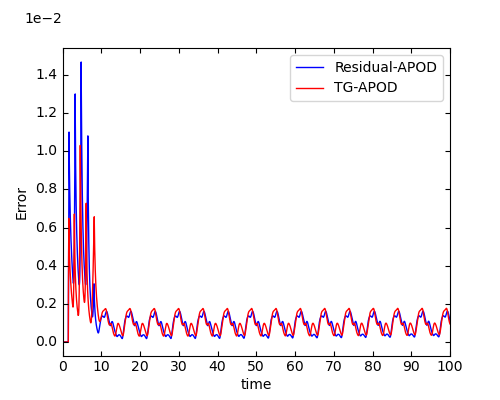}
		\end{minipage}%
	}%
	\vspace{-1em}
	\caption*{(a).$\epsilon=0.5$}
	\vspace{-1em}
	\subfigure{
		\begin{minipage}[t]{0.45\linewidth}
			\centering
			\includegraphics[width=6cm,height=4cm]{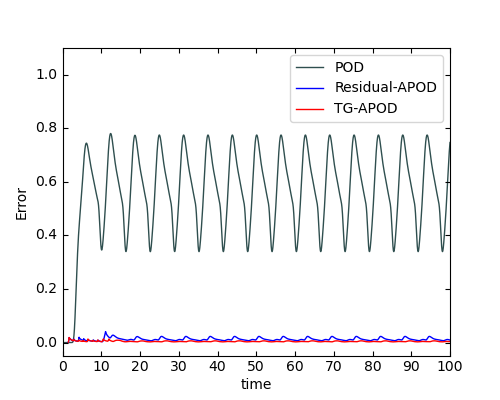}
		\end{minipage}
	}%
	\subfigure{
		\begin{minipage}[t]{0.45\linewidth}
			\centering
			\includegraphics[width=6cm,height=4cm]{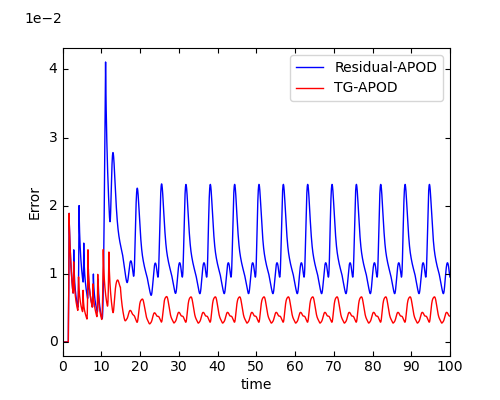}
		\end{minipage}
	}%
	\vspace{-1em}
	\caption*{(b).$\epsilon=0.1$}
	\vspace{-1em}
	\subfigure{
		\begin{minipage}[t]{0.45\linewidth}
			\centering
			\includegraphics[width=6cm,height=4cm]{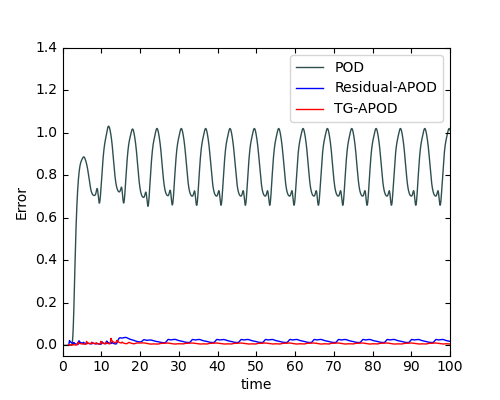}	
		\end{minipage}%
	}
	\subfigure{
		\begin{minipage}[t]{0.45\linewidth}
			\centering
			\includegraphics[width=6cm,height=4cm]{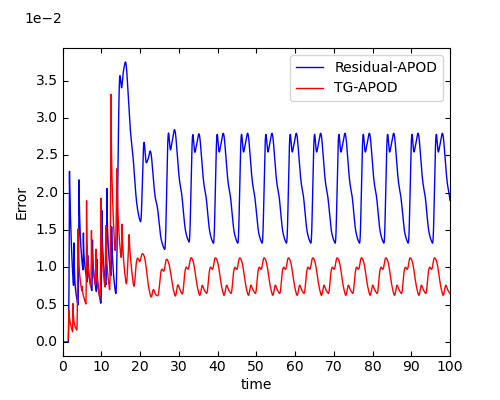}
		\end{minipage}%
	}%
	\vspace{-1em}
	\caption*{(c).$\epsilon=0.05$}	
	\vspace{-1em}
	\subfigure{
		\begin{minipage}[t]{0.45\linewidth}
			\centering
			\includegraphics[width=6cm,height=4cm]{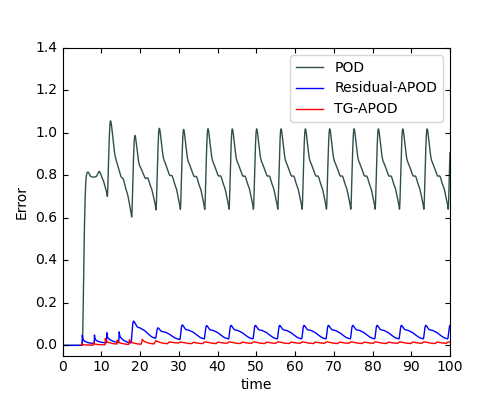}	
		\end{minipage}%
	}%
	\subfigure{
		\begin{minipage}[t]{0.45\linewidth}
			\centering
			\includegraphics[width=6cm,height=4cm]{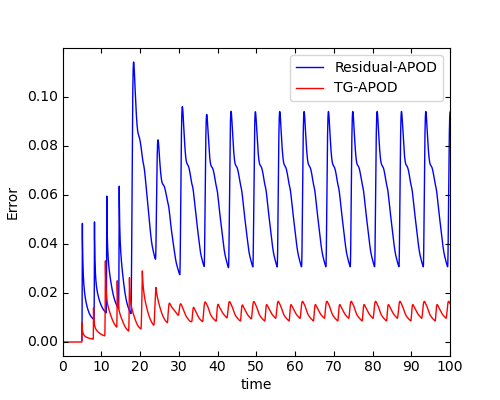}
		\end{minipage}%
	}%
	\vspace{-1em}
	\caption*{(d).$\epsilon=0.01$}	
	\vspace{-1em}
	\centering
	\caption{The  error curves for solution of (\ref{Kolmogorov-equation}) with different $\epsilon$ obtained by the standard POD method, the  residual based adaptive POD method, and the two-grid based adaptive POD method, respectively. }	
	\label{compare-twomethod-figure1}
\end{figure}

To see more clearly, we compare the error of the numerical solutions obtained by the different methods in Fig. \ref{compare-twomethod-figure1}. 

 We should point out that among the different results obtained by residual-APOD with different $\eta_0$, we choose the best one to show in Fig.\ref{compare-twomethod-figure2}. 

In Fig. \ref{compare-twomethod-figure1}, the x-axis is time, the y-axis is relative error of numerical solution. The results obtained by the standard POD method, the residual based adaptive POD algorithm, and our two-grid based adaptive POD algorithm are reported in lines with color darkslategray, blue,  and red, respectively.

From the left figure of  Fig. \ref{compare-twomethod-figure1}, we can see that the  error curve  obtained by the standard POD method is above both those obtained by the two adaptive POD methods, which means that the adaptive POD methods are more efficient than the standard POD method. From the right figure of Fig.2, we can see that our two-grid based adaptive POD method is more efficient than the residual based adaptive POD method.

\subsection{Arnold-Beltrami-Childress (ABC) flow} We then consider the ABC flow. The ABC flow was introduced by Arnold, Beltrami, and Childress \cite{Kolmogorov-flow} to study chaotic advection, enhanced transport and dynamo effect,
 see \cite{biferale1995eddy,galanti1992linear,wirth1995eddy,xin2016periodic} for the details. 
 
 We consider the following equation
\begin{eqnarray}\label{ABC-equation}
\left\{
\begin{aligned}
&u_{t}-\epsilon \Delta u +\vec B(x,y,z,t) \cdot \nabla u = f(x,y,z,t), \quad (x, y, z)\in\Omega , t\in [0, T],\\
&u(x,y,z,0)= 0, \\
&u(x+2\pi,y,z,t)=u(x,y+2\pi,z,t)=u(x,y,z+2\pi,t)=u(x,y,z,t),\\
\end{aligned}
\right.
\end{eqnarray}
where
\begin{eqnarray*}
	\vec B(x,y,z,t)&=&( \sin (z +  \sin w t )  + \cos ( y +  \sin w t ),  \sin( x +  \sin w t) \\
	&&+ \cos (z +  \sin w t), \sin (y +  \sin w t) + \cos (x + \sin w t) ), \\
	f(x,y,z,t)&=&-\sin (z +  \sin w t )  - \cos ( y +  \sin w t ),\\
	\Omega = [0, 2\pi&]^3 &,  T=100.
\end{eqnarray*}
 
For this example, we also test $4$ different cases with  $\epsilon=0.5, 0.1$, $0.05$, and $0.01$, respectively.
We divide the domain $[0, 2\pi]^3$ into tetrahedrons to get the initial grid containing $6$ elements, then refine the initial mesh $22$ times uniformly using bisection to get the fine mesh.
We set $w=1.0$, $\delta t =0.005$, and choose the finite element basis to be piecewise linear function.
 For the two cases with  $\epsilon=0.5, 0.1$, the parameters are chosen as $T_0=1.5, \delta T=1, \delta M=5$; for other cases with $\epsilon=0.05,0.01$, the parameters are chosen as $T_0= 5, \delta T=4,\delta M=20$.
In the two-grid based adaptive POD algorithm, we refine the initial mesh $16$ times uniformly using bisection to obtain the coarse mesh, and the time step of the coarse mesh is set to $0.09$. We use $36$ processors for the simulation.

We use the similar choice for parameters $\gamma_i$ ($i=1,2,3$) and $\eta_0$ as for the  Kolmogorov flow. That is,  for our  two-grid based adaptive POD method, we choose them as $\gamma_1=\gamma_2 = 0.999$, $\gamma_3 = 1.0-1.0\times 10^{-8}$, and $\eta_0 = 0.005$, while   for the residual-APOD method, we choose  $\gamma_1=\gamma_2 = 0.999$, $\gamma_3 = 1.0-1.0\times 10^{-8}$, but setting  different values for $\eta_0$ to find out the best performance of the residual-APOD method.

\begin{figure}[!htbp]	
		\begin{minipage}[t]{0.5\textwidth}
			\centering
			\includegraphics[width =\textwidth]{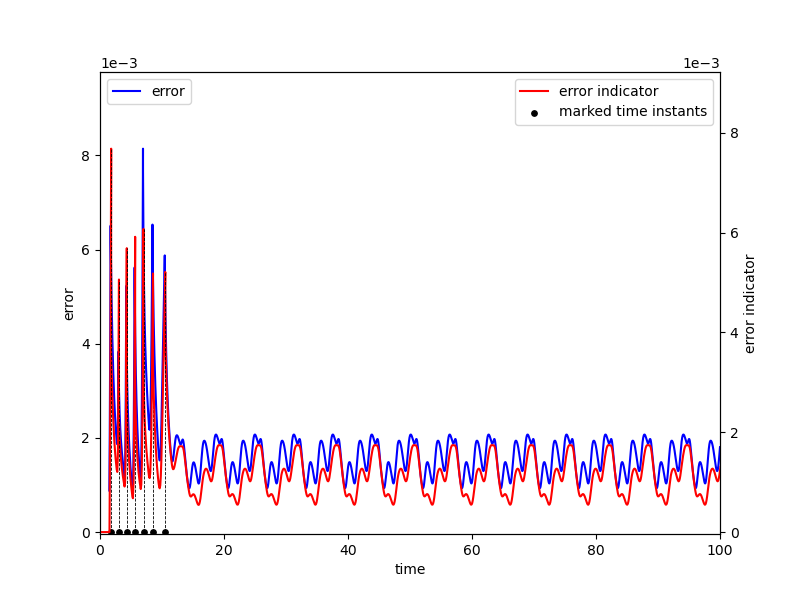}
			\caption*{a).$\epsilon=0.5$.}
		\end{minipage}
		\begin{minipage}[t]{0.5\textwidth}
			\centering
			\includegraphics[width = \textwidth]{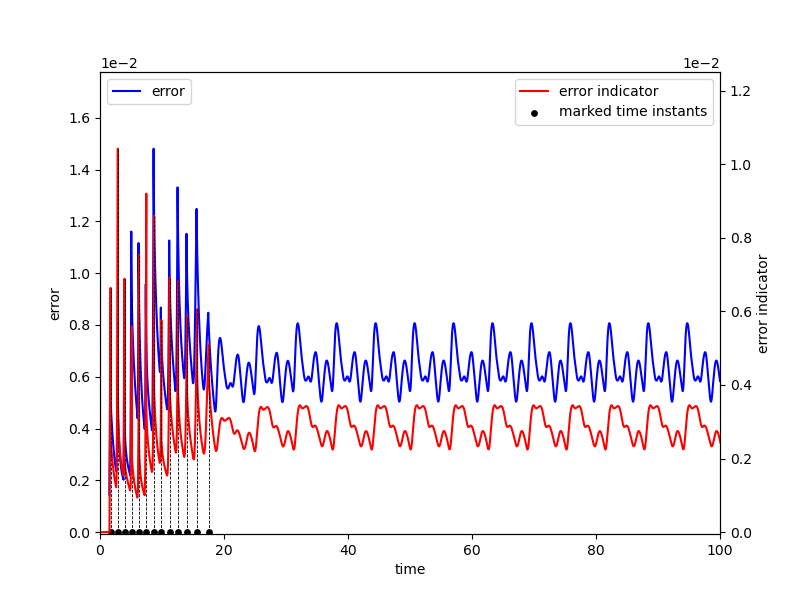}
			\caption*{b).$\epsilon=0.1$.}
		\end{minipage}
		\begin{minipage}[t]{0.5\textwidth}
			\centering
			\includegraphics[width = \textwidth]{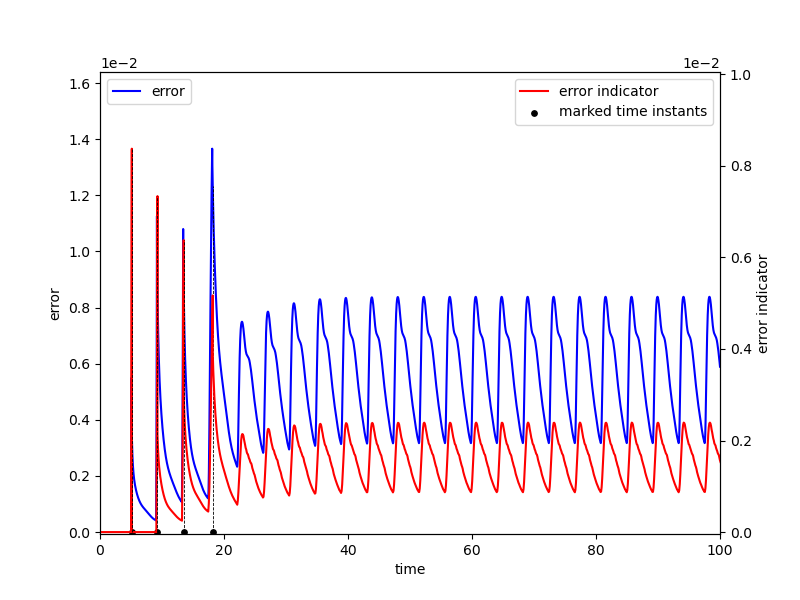}
			\caption*{c).$\epsilon=0.05$.}
		\end{minipage}	
		\begin{minipage}[t]{0.5\textwidth}
			\centering
			\includegraphics[width =\textwidth]{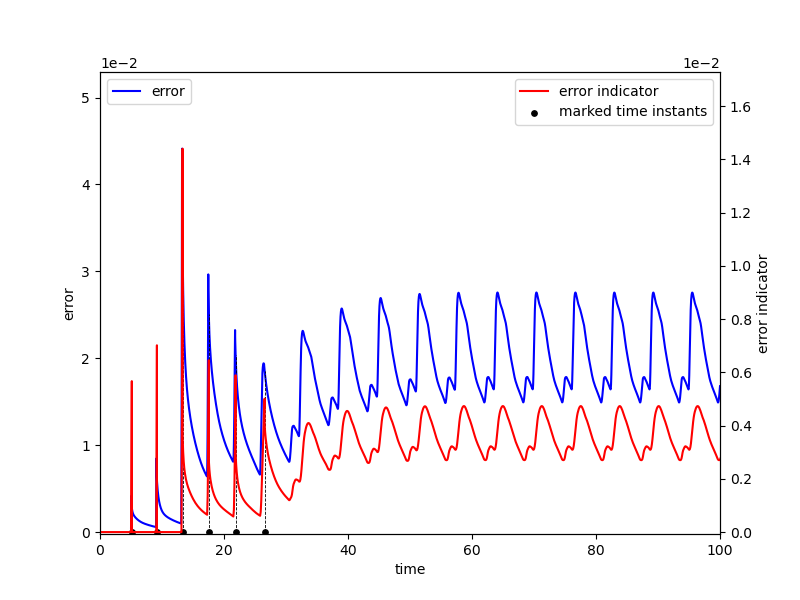}
			\caption*{d).$\epsilon=0.01$.}
		\end{minipage}
		\caption{The evolution curves  of the error indicator and the error with $\epsilon=0.5, 0.1, 0.05,0.01$.}
		\label{err-comp-2}
	\end{figure}

	Similar to the Kolmogorov flow, we first use some numerical results to show the efficiency of our error indicator and the mark startegy. Fig \ref{err-comp-2} 
	shows the evolution curves of the error indicator and the error, and the time instance being marked.    
	
	Similarly,  by  comparing the  curves of the error indicator and the error, and the  curve of the error and the time instance being marked to be updated in Fig. \ref{err-comp-2}, we know that  the error indicator is efficient.  
	
	We then compare the numerical results  obtained by our two-grid based adaptive POD method, the residual-APOD method, the standard POD method. The detailed results are listed in Table \ref{compare3method-table2}.  The results obtained by standard finite element method are also provided as a reference. The notation in Table 2 has the same meaning as in Table \ref{compare3method-table1}.

\begin{table}[!htbp]
	\centering
	\caption{The results of the ABC flow  with different $\epsilon$ obtained by Standard-FEM, POD, Residual-APOD, and TG-APOD method, respectively.}
	\begin{tabular}{|c|c| c |c| c | c | c |}
		\hlinew{1.5pt}
		$\epsilon $       &   Methods & $\eta_0$ & DOFs& Average Error & Times (s) \\
		\hlinew{1.5pt}
		\multirow{5}{*}{0.5}  & Standard FEM & -- &4194304 & -- & 16690.55\\
		&POD   & --  &  12 &  0.303585      &  534.96\\
		&Res-APOD  &$5\times 10^{-3}$ & 9  & 0.191637  &  7749.15\\
		&Res-APOD  &$1\times 10^{-3}$ & 25 & 0.007640 & 8636.51 \\
		&TG-APOD   &$5\times 10^{-3}$ &  33  &   0.001660    &2240.67 \\
		\hlinew{1.5pt}
		\multirow{5}{*}{0.1}  & FEM  &-- & 4194304 &-- & 15238.87 \\
		&POD   & --  &14  &  0.658382  & 536.73 \\
		&Res-APOD  &$5\times 10^{-3}$ & 32   & 0.196853 & 9041.01\\
		&Res-APOD  & $1\times 10^{-3}$ & 57  & 0.024187 & 10213.65 \\
		&TG-APOD  &$5\times 10^{-3}$ & 74  & 0.006207 & 4406.87\\
		\hlinew{1.5pt}
		\multirow{5}{*}{0.05}    & FEM   &--& 4194304 & -- & 14439.24\\
		&POD   &-- &   42 &    0.504838     &1485.62\\
		&Res-APOD &$5\times 10^{-3}$ & 54 & 0.059658 &  10094.97\\
		&Res-APOD  &$1\times 10^{-3}$ & 67 & 0.014872 & 10734.27 \\
		&TG-APOD &$5\times 10^{-3}$  &  100 &  0.005077  & 5011.81 \\
		\hlinew{1.5pt}
		\multirow{5}{*}{0.01}    &  FEM  &-- &4194304 &-- & 16074.96\\
		&POD   &-- & 61 &    0.676701& 1878.01 \\
		&Res-APOD  & $1\times 10^{-2}$   & 64 & 0.172755 & 10275.33 \\
		&Res-APOD &$5\times 10^{-3}$  &140 & 0.042498 & 13322.86  \\
		&TG-APOD  &$5\times 10^{-3}$ & 192   & 0.015851  & 10279.38\\ \hlinew{1.5pt}
	\end{tabular}\label{compare3method-table2}
\end{table}

\indent Similar to the first example, we can see from Table \ref{compare3method-table2} that those POD type methods  can indeed reduce the number of basis a lot. For this example, the accuracy for the approximation obtained by the standard POD method is so large that the results are meaningless.  While  both the residual based adaptive POD method and our two-grid based adaptive POD method can produce numerical solutions with high accuracy. For  the two adaptive POD methods,  our two-grid adaptive POD method can yield  higher accuracy numerical solutions than the residual based adaptive POD method. When it comes to cpu time,
our two-grid adaptive method takes less cpu time than the residual based adaptive POD method. These show that our two-grid adaptive POD method is more efficient.

\begin{figure}[!htbp]
	\centering	
	\subfigure{
		\begin{minipage}[t]{0.45\linewidth}
			\centering
			\includegraphics[width=6cm,height=4cm]{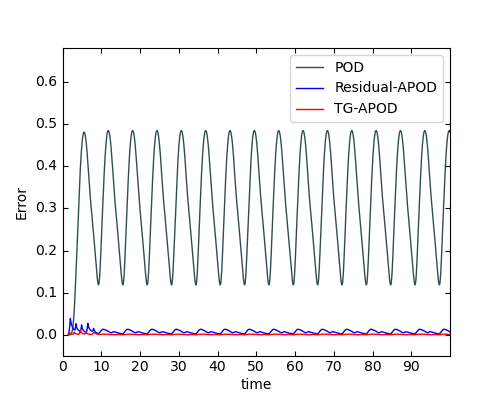}	
		\end{minipage}%
	}%
	\subfigure{
		\begin{minipage}[t]{0.45\linewidth}
			\centering
			\includegraphics[width=6cm,height=4cm]{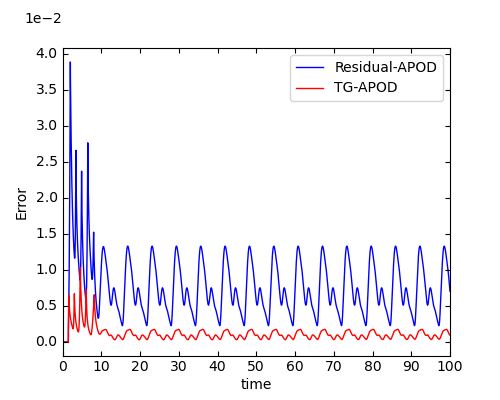}
		\end{minipage}%
	}%
	\vspace{-1em}
	\caption*{(a).$\epsilon=0.5$}
	\vspace{-1em}
	\subfigure{
		\begin{minipage}[t]{0.45\linewidth}
			\centering
			\includegraphics[width=6cm,height=4cm]{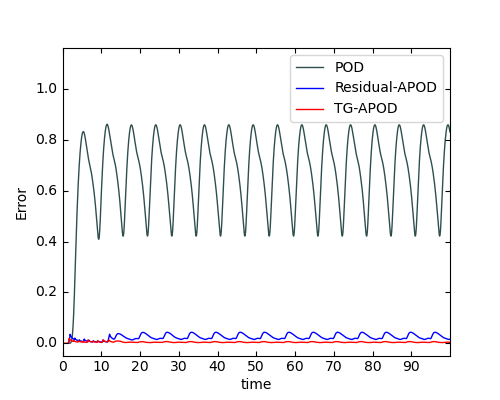}
		\end{minipage}
	}%
	\subfigure{
		\begin{minipage}[t]{0.45\linewidth}
			\centering
			\includegraphics[width=6cm,height=4cm]{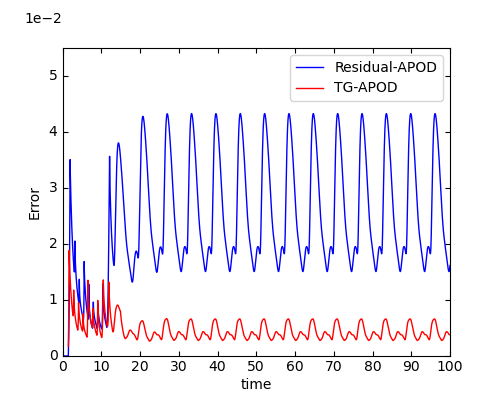}
		\end{minipage}
	}%
	\vspace{-1em}
	\caption*{(b).$\epsilon=0.1$}
	\vspace{-1em}
	\subfigure{
		\begin{minipage}[t]{0.45\linewidth}
			\centering
			\includegraphics[width=6cm,height=4cm]{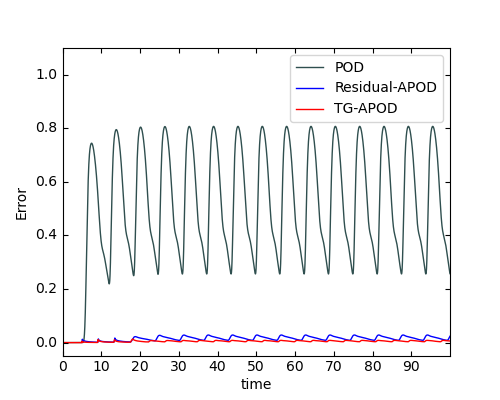}	
		\end{minipage}%
	}%
	\subfigure{
		\begin{minipage}[t]{0.45\linewidth}
			\centering
			\includegraphics[width=6cm,height=4cm]{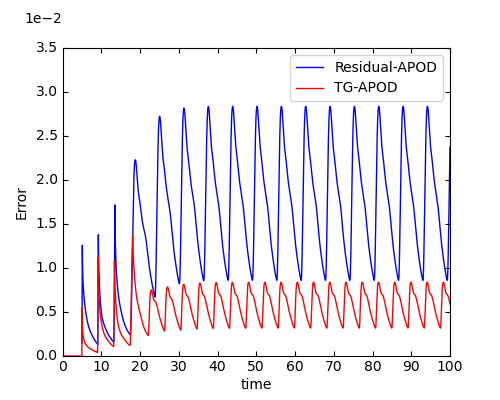}
		\end{minipage}%
	}%
	\vspace{-1em}
	\caption*{(c).$\epsilon=0.05$}	
	\vspace{-1em}
	\subfigure{
		\begin{minipage}[t]{0.45\linewidth}
			\centering
			\includegraphics[width=6cm,height=4cm]{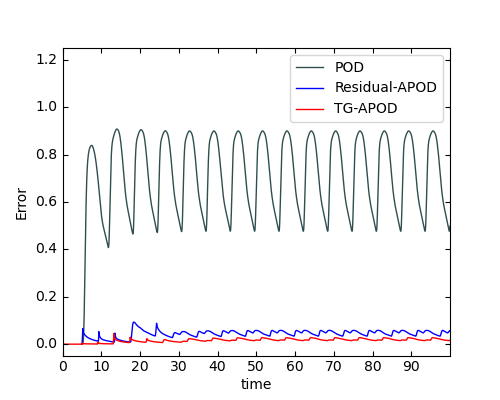}	
		\end{minipage}%
	}%
	\subfigure{
		\begin{minipage}[t]{0.45\linewidth}
			\centering
			\includegraphics[width=6cm,height=4cm]{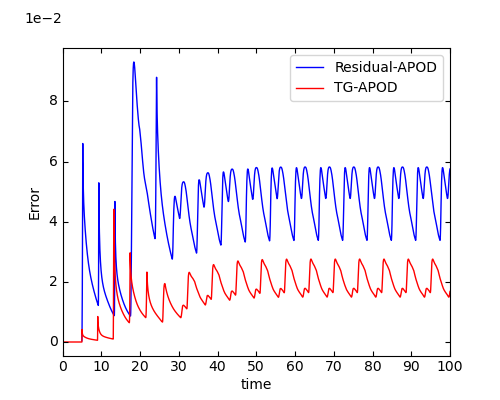}
		\end{minipage}%
	}%
	\vspace{-1em}
	\caption*{(d).$\epsilon=0.01$}	
	\vspace{-1em}
	\centering
	\caption{The evolution curves of error for solution of (\ref{ABC-equation}) with different $\epsilon$ obtained by the standard POD method, the  residual based adaptive POD method, and the two-grid based adaptive POD method, respectively. }	
	\label{compare-twomethod-figure2}
\end{figure}

Similarly, we  show the error curves of the numerical solutions obtained by the different methods in Fig. \ref{compare-twomethod-figure2}.

In Fig. \ref{compare-twomethod-figure2}. the x-axis is time, the y-axis is the relative error of the POD approximations. The results obtained by the standard POD method, the residual based adaptive POD method, and our two-grid based adaptive POD method are reported in line with color darkslategray, blue,  and red, respectively.

From Fig. \ref{compare-twomethod-figure2}, we can see more clearly that both the two-grid adaptive POD method and the residual based adaptive POD method behave much better than the standard POD method,
 and our two-grid based adaptive POD method behaves a little better than the residual based adaptive POD method.

  We have done some more  numerical experiments for setting different parameters $w$. The results for cases with $w= 1.5, 2 ,2.5$ are  shown in Table \ref{D_0&w}.

\begin{table}[!htbp]
	\centering
	\small
	\caption{ Results of the ABC flow  with different $w$ and $\epsilon$ obtained by standard-FEM,POD,  Residual-APOD and TG-APOD algorithm, respectively}
	\label{D_0&w}
	\begin{tabular}{|c|c| c |c| c|c | c | c |}
		\hlinew{1.5pt}
		$\epsilon$       & $w$   &  Methods & $\eta_0$& DOFs& Average Error  & Time (s) \\
		\hline
		\multirow{5}{*}{0.5}  & \multirow{5}{*}{1.5} &FEM  & --&4194304 &--& 16745.85\\
		&                     &POD &--  & 14 &  0.294394  & 561.23\\
		&                     &Res-APOD  & $5\times 10^{-3}$ &  10 &  0.138573  &  7773.07\\
		&                           &Res-APOD  & $5\times 10^{-4}$ & 30 & 0.003371 & 9017.68  \\
		&                     &TG-APOD  & $5\times 10^{-3}$ & 30 &  0.002611    & 2184.55 \\
		\hline
		\multirow{5}{*}{0.5}    & \multirow{5}{*}{2.0}  & FEM&--& 4194304 &--&17258.89\\
		&                             &POD   & -- &  15  &  0.270138    & 530.91 \\
		&                           &Res-APOD & $5\times 10^{-3}$ & 6  &  0.336776  & 7581.42 \\
		&                           &Res-APOD  & $1\times 10^{-3}$ & 21 &  0.013750 &  8497.94 \\
		&                           &TG-APOD  & $5\times 10^{-3}$& 31 &  0.001515  & 1927.63 \\
		\hline
		\multirow{5}{*}{0.5}    & \multirow{5}{*}{2.5}& FEM & -- &4194304& --& 17161.84 \\
		&                             &POD  & --&  17&    0.216109    &  602.85 \\
		&                           &Res-APOD & $5\times 10^{-3}$ &   6  & 0.338068 & 7745.82 \\
		&                           &Res-APOD  & $5\times 10^{-4}$ & 21 & 0.008787 &  8440.45 \\
		&                           &TG-APOD  & $5\times 10^{-3}$ &   26   & 0.002678 & 1617.53\\
		\hlinew{1.5pt}
		\multirow{5}{*}{0.1}    & \multirow{5}{*}{1.5} & FEM &--&4194304&--& 15232.78 \\
		&                             &POD  & --&15 &   0.729666  &  531.42  \\
		&                           &Res-APOD  &$5\times 10^{-3}$ & 32  &  0.148178  &  8576.58 \\
		&                           &Res-APOD  & $5\times 10^{-4}$ & 65 &0.012987  &  10359.17 \\
		&                           &TG-APOD   & $5\times 10^{-3}$& 78 & 0.006964 & 4238.26\\
		\hline
		\multirow{5}{*}{0.1}    & \multirow{5}{*}{2.0} & FEM  & --&4194304& --& 15324.38 \\
		&                             &POD & -- &16 & 0.751065 & 536.24  \\
		&                           &Res-APOD & $5\times 10^{-3}$& 27  &0.161848 & 8339.32\\
		&                           &Res-APOD  & $5\times 10^{-4}$  & 57 & 0.018898 &  9790.15 \\
		&                           &TG-APOD  & $5\times 10^{-3}$ & 68   & 0.006915  & 3710.35\\
		\hline
		\multirow{5}{*}{0.1}    & \multirow{5}{*}{2.5}         & FEM   & -- & 4194304  &    --     &  15104.35 \\
		&                             &POD  & --& 18 &0.657750   & 589.66\\
		&        &Res-APOD   &$5\times 10^{-3}$  & 21 & 0.222043  & 7886.83\\
		&                           &Res-APOD  &$1\times 10^{-3}$  & 40 & 0.041986 & 9071.88  \\
		&        &TG-APOD& $5\times 10^{-3}$  & 75  &  0.005646  & 4212.24 \\
		\hlinew{1.5pt}
		\multirow{5}{*}{0.05}    & \multirow{5}{*}{1.5}         & FEM  & --& 4194304   & -- &14141.53 \\
		&                             &POD&--&45 &0.445832&1460.26 \\
		&                      &Res-APOD  &  $5\times 10^{-3}$    &  42 &   0.148505  &  9410.69 \\
		&                           &Res-APOD  & $2\times 10^{-3}$ & 62 &  0.041819 & 10338.44  \\
		&                           &TG-APOD & $5\times 10^{-3}$ &  98 &   0.004177 & 5361.94 \\
		\hline
		\multirow{5}{*}{0.05}    & \multirow{5}{*}{2.0}    & FEM  &-- & 4194304 & -- &14342.82 \\
		&                             &POD &--  & 47 &  0.325020  &1425.27 \\
		&                           &Res-APOD & $5\times 10^{-3}$ & 44 &  0.088758 &  9303.73\\
		&                           &Res-APOD  & $1\times 10^{-3}$ & 63  & 0.022060 &  10014.03 \\
		&                           &TG-APOD & $5\times 10^{-3}$  & 83 &0.004348  & 4299.98 \\
		
		\hline
		\multirow{5}{*}{0.05}    & \multirow{5}{*}{2.5} & FEM & --&4194304&--& 14235.11\\
		&                             &POD  &-- & 49 &     0.245578   & 1528.92 \\
		&                           &Res-APOD & $5\times 10^{-3}$ &  46  & 0.057361 &  9220.55\\
		&                           &Res-APOD  & $5\times 10^{-4}$ & 63 &0.016595  &  10375.65 \\
		&                           &TG-APOD &$5\times 10^{-3}$  & 84 & 0.007256 &  4313.27\\
		\hlinew{1.5pt}
		\multirow{5}{*}{0.01}    & \multirow{5}{*}{1.5}         & FEM &-- &4194304  &--&  16223.32 \\
		&                             &POD  &-- &46&0.670506 & 1533.39 \\
		&                           &Res-APOD  &  $2\times 10^{-2}$  & 95 & 0.155503  &  11154.77 \\
		&                           &Res-APOD &  $5\times 10^{-3}$ & 123 & 0.060578 &  12527.66  \\
		&                           &TG-APOD & $5\times 10^{-3}$  & 187 & 0.016889  & 9337.59\\
		\hline
		\multirow{5}{*}{0.01}    & \multirow{5}{*}{2.0}         & FEM  &-- &4194304&--&16376.39\\
		&                             &POD   &--& 48 & 0.557215  &1638.90 \\
		&                           &Res-APOD  & $8\times 10^{-3}$ & 97 & 0.122302 &  10954.07  \\&                           &Res-APOD & $5\times 10^{-3}$ & 130 &  0.047402     &  12396.49 \\
		&                           &TG-APOD  & $5\times 10^{-3}$ &194 &   0.010994  & 9371.19 \\
		\hline
		\multirow{5}{*}{0.01}    & \multirow{5}{*}{2.5}         & FEM   &--& 4194303 & --&   16226.39\\
		&                             &POD&-- & 49& 0.582593   &1589.38  \\
		&                           &Res-APOD  & $1\times 10^{-2}$ & 64 & 0.172755 &  9645.38 \\
		&                           &Res-APOD  & $5\times 10^{-3}$ & 97 & 0.080639 &  11499.29 \\
				&                           &TG-APOD& $5\times 10^{-3}$ &188  &0.015246  & 9427.31\\
		\hlinew{1.5pt}
	\end{tabular}
\end{table}

From Table \ref{D_0&w}, we can obtain the same conclusion as those from Table \ref{compare3method-table1} and Table \ref{compare3method-table2}, that is, adaptive POD methods outperform the POD method a lot, and our two-grid based adaptive method outperforms the residual based adaptive method.

\section{Concluding remarks}\label{sec-conclude}
In this paper, we proposed a two-grid based adaptive POD method to solve the time dependent partial differential equations. We apply our method to some typical 3D advection-diffusion equations,  the Kolmogorov flow and the ABC flow. Numerical results show that our two-grid based adaptive POD method is more efficient than the residual based  adaptive POD method, especially is much more efficient than the standard POD method. We should also point out that for problems with extreme sharp phase change,  our method may behave not so well.  Here, we simply use the relative error of the POD solution on the coarse spacial and temporal mesh to construct the error indicator. In our future work, we will consider  other error indicator based on our two-grid approach or some other approach, and study other types of time dependent partial differential equations. 

\begin{acknowledgements}
	The authors would like to thank the anonymous referees for their nice and useful comments and suggestions
	that improve the quality of this paper.
\end{acknowledgements}

\begin{appendix}
\renewcommand{\appendixname}{Appendix~\Alph{section}}

\section{Numerical experiments for tuning parameters $\gamma_i(i=1,2,3)$}\label{numer-tune-gamma}

	In this section, we do some numerical experiments for tuning the parameters $\gamma_i(i=1,2,3)$ to find out a good choice. 
	
	Since the parameters $\gamma_1$ and $\gamma_2$ are both for extracting POD modes from the snapshots obtained from the standard finite element approximation of the original dynamic system, we set them to be the same value, that is, $\gamma_1 = \gamma_2$. We have tested the Kolmogorov flow with  $\epsilon=0.05$ and $\epsilon=0.01$ respectively, and the ABC flow with $\epsilon=0.01$ by using the following different parameters:
	
	\begin{itemize}
		\item[1)] case 1: $\gamma_1= \gamma_2 = 1.0-1.0\times 10^{-8}$, set $\gamma_3$ to be $0.9, 0.99, 0.999,$ and $0.9999$ respectively.  
		\item[2)] case 2: $\gamma_3 = 1.0-1.0\times 10^{-8}$, set $\gamma_1$ and $\gamma_2$ to be  $0.9, 0.99, 0.999,$ and $0.9999$ respectively.
		\item[3)] case 3: $\gamma_1=\gamma_2 = \gamma_3 = 1.0-1.0\times 10^{-2}, 1.0-1.0\times 10^{-4}, 1.0-1.0\times 10^{-6}, 1.0-1.0\times 10^{-8}$, respectively. 
	\end{itemize}

	In our experiments, we choose $\eta_0 = 0.005$. The other parameters are same as stated in Section 4.

	{\bf   Kolmogorov flow with $\epsilon = 0.01$ }
	
	We first see the results for the Kolmogorov flow with $\epsilon = 0.01$.  
	
	The results for  case 1 are shown in  
	Table \ref{Tab-Kolmogrov-0.01-gamma3}. We can see that for fixing $\gamma_1$ and $\gamma_2$, the smaller  $\gamma_3$, the more accurate the results, which means $\gamma_3$ close to $1$ may be a good choice. Then,  we can see that when $\gamma_1$ and $\gamma_2$ are chosen to be too close to $1$, no matter how to choose the parameter $\gamma_3$, the cpu-time cost is too much. These results tell us that setting  $\gamma_1$ and $\gamma_2$  too close to $1$ is not a good choice.

	\begin{table}[!htbp]
		\centering
		\caption{Detailed information for results of the Kolmogorov flow with $\epsilon=0.01$ obtained by setting $\gamma_1=\gamma_2=1.0-1.0\times 10^{-8}$ and  $\gamma_3=0.9,0.99, 0.999, 0.9999$, respectively.}
		\label{Tab-Kolmogrov-0.01-gamma3}
		\begin{tabular}{| c| c|c | c | c |c|}
			\hlinew{1.2pt}
			&POD-updates&   DOFs   &   Average Error&  Times (s) \\
			\hlinew{1.2pt}
			$\gamma_3=0.9$   & 30 & 66 &0.059832  & 15770.89 \\
			$\gamma_3=0.99$    & 30 & 236 & 0.027262 & 31568.01  \\
			$\gamma_3=0.999$   & 22 &339 & 0.018255  & 38263.16 \\
			$\gamma_3=0.9999$ & 6 & 222 & 0.043950 &   9501.85 \\
			\hlinew{1.5pt}
		\end{tabular}
	\end{table}

	The results for  case 2 are shown in  
	Table \ref{Tab-Kolmogrov-0.01-gamma1-2}, from which  we can see that for fixed $\gamma_3$, which is very close to $1$, setting $\gamma_1$ and $\gamma_2$ to be $0.999$  can obtain results with high accuracy. Taking into account the cpu time cost,  $0.999$ is a better choice.

	\begin{table}[!htbp]
		\centering
		\caption{Detailed information for results of the Kolmogorov flow with $\epsilon=0.01$ obtained by setting $\gamma_3=1.0-1.0\times 10^{-8}$ and   $\gamma_1=\gamma_2=0.9,0.99, 0.999, 0.9999$, respectively.}
		\label{Tab-Kolmogrov-0.01-gamma1-2}
		\begin{tabular}{| c| c|c | c | c |c|}
			\hlinew{1.2pt}
			&POD-updates&   DOFs   &   Average Error&  Times (s) \\
			\hlinew{1.2pt}
			$\gamma_1=\gamma_2=0.9$   & 17 & 124&0.031398 & 7140.06\\
			$\gamma_1=\gamma_2=0.99$    & 7 & 108 & 0.014088 & 4322.85 \\
			$\gamma_1=\gamma_2=0.999$   & 7 &138 &0.010969  &  6299.70   \\
			$\gamma_1=\gamma_2=0.9999$ & 7 & 185 & 0.015274   & 8457.13   \\
			\hlinew{1.5pt}
		\end{tabular}
	\end{table}

	The results for  case 3 are shown in  
	Table \ref{Tab-Kolmogrov-0.01-gamma}.     
	We see that taking into account both the accuracy and the cpu time cost, these parameters are not as good as those for case 2.

	\begin{table}[!htbp]
		\centering
		\caption{Detailed information for results of the Kolmogorov flow with $\epsilon=0.01$ obtained by setting  $\gamma_1=\gamma_2=\gamma_3=\gamma=1.0-1.0\times 10^{-2},1.0-1.0\times 10^{-4},1.0-1.0\times 10^{-6},1.0-1.0\times 10^{-8}$ respectively.}
		\label{Tab-Kolmogrov-0.01-gamma}
		\begin{tabular}{| c| c |c| c | c |c|}
			\hlinew{1.2pt}
			&  POD-updates	&POD-DOFs&  Average Error& Times (s) \\
			\hlinew{1.2pt}
			$\gamma=1-10^{-2}$   &31 &55 & 0.013877  & 15354.34\\
			$\gamma=1-10^{-4}$    & 7 & 182   &  0.012155   &  8351.72 \\
			$\gamma=1-10^{-6}$  & 6 &204  & 0.024557 & 8826.03 \\
			$\gamma=1-10^{-8}$ & 6 & 223 & 0.043659  &  9585.60 \\
			\hlinew{1.5pt}
		\end{tabular}
	\end{table}

	Anyway, by comparing Table \ref{Tab-Kolmogrov-0.01-gamma3}, Table \ref{Tab-Kolmogrov-0.01-gamma1-2}, Table \ref{Tab-Kolmogrov-0.01-gamma}, we can see that the case of $\gamma_1=\gamma_2=0.999$, $\gamma_3 = 1.0-1.0\times 10^{-8}$ is the best choice among all these cases, if taking both the accuracy and the cpu-time cost into account. 
	 
	{\bf  the Kolmogorov flow with $\epsilon = 0.05$ }
	
	We then take a look at the results for the Kolmogorov flow with $\epsilon = 0.05$, which are shown in  Table\ref{Tab-Kolmogrov-0.05-gamma3}, Table \ref{Tab-Kolmogrov-0.05-gamma1-2}, and
	Table \ref{Tab-Kolmogrov-0.05-gamma}.

	\begin{table}[!htbp]
		\centering
		\caption{Detailed information for results of the Kolmogorov flow with $\epsilon=0.05$ obtained by setting $\gamma_1=\gamma_2=1.0-1.0\times 10^{-8}$ and  $\gamma_3=0.9,0.99, 0.999, 0.9999$, respectively.}
		\label{Tab-Kolmogrov-0.05-gamma3}
		\begin{tabular}{| c| c|c | c | c |c|}
			\hlinew{1.2pt}
			&POD-updates&   DOFs   &   Average Error&  Times (s) \\
			\hlinew{1.2pt}
		$\gamma_3=0.9$   & 92 & 86 & 0.159828 & 16482.18\\
		$\gamma_3=0.99$    & 92  & 473 & 0.076002 & 52158.13\\
		$\gamma_3=0.999$   & 65 & 658 & 0.031944  & 86469.45 \\
		$\gamma_3=0.9999$ & 10 & 182 & 0.007093 &   7415.55 \\
			\hlinew{1.5pt}
		\end{tabular}
	\end{table}

	\begin{table}[!htbp]
		\centering
		\caption{Detailed information for results of the Kolmogorov flow with $\epsilon=0.05$ obtained by setting $\gamma_3=1.0-1.0\times 10^{-8}$ and   $\gamma_1=\gamma_2=0.9,0.99, 0.999, 0.9999$, respectively.}
		\label{Tab-Kolmogrov-0.05-gamma1-2}
		\begin{tabular}{| c| c|c | c | c |c|}
			\hlinew{1.2pt}
			&POD-updates&   DOFs   &   Average Error&  Times (s) \\
			\hlinew{1.2pt}
		$\gamma_1=\gamma_2=0.9$   & 17 & 37& 0.028841 & 3114.22\\
	$\gamma_1=\gamma_2=0.99$    & 14 & 59 & 0.012810 & 3645.05 \\
	$\gamma_1=\gamma_2=0.999$   &13  & 83& 0.008663 &  3959.62 \\
	$\gamma_1=\gamma_2=0.9999$ & 12 & 109 &0.007722  &  4690.68  \\
		\hlinew{1.5pt}
		\end{tabular}
	\end{table}

	\begin{table}[!htbp]
		\centering
		\caption{Detailed information for results of the Kolmogorov flow with $\epsilon=0.05$ obtained by setting  $\gamma_1=\gamma_2=\gamma_3=\gamma=1.0-1.0\times 10^{-2},1.0-1.0\times 10^{-4},1.0-1.0\times 10^{-6},1.0-1.0\times 10^{-8}$ respectively.}
		\label{Tab-Kolmogrov-0.05-gamma}
		\begin{tabular}{| c| c |c| c | c |c|}
			\hlinew{1.2pt}
			&  POD-updates	&POD-DOFs&  Average Error& Times (s) \\
			\hlinew{1.2pt}
			$\gamma=1-10^{-2}$   & 92 & 49 &  0.023817 & 12711.97\\
		$\gamma=1-10^{-4}$    & 13 & 104 & 0.006539 & 4726.16  \\
		$\gamma=1-10^{-6}$  & 11 &152  & 0.006474 & 6421.31  \\
		$\gamma=1-10^{-8}$ &10  & 188 & 0.007049  &   7739.79 \\
				\hlinew{1.5pt}
		\end{tabular}
	\end{table}

	By comparing Table \ref{Tab-Kolmogrov-0.05-gamma3}, Table \ref{Tab-Kolmogrov-0.05-gamma1-2}, Table \ref{Tab-Kolmogrov-0.05-gamma}, we can also see that the case of $\gamma_1=\gamma_2=0.999$, $\gamma_3 = 1.0-1.0\times 10^{-8}$ is the best choice among all these cases if taking both the accuracy and the cpu-time cost into account.

{\bf ABC flow with $\epsilon = 0.01$}

We have also done some numerical experiments for the ABC flow with $\epsilon = 0.01$. For this example, the results obtained by setting different choice of the parameters $\gamma_i(i=1,2,3)$ are shown in Table \ref{Tab-ABC-0.01-gamma3}, Table \ref{Tab-ABC-0.01-gamma1-2}, and Table \ref{Tab-ABC-0.01-gamma}. From these results, we can  obtain the same conclusion as those obtained from the Kolmogorov flow, that is,  $\gamma_1=\gamma_2=0.999$, $\gamma_3 = 1.0-1.0\times 10^{-8}$ is the best choice among all these cases if taking into account both the accuracy and the cpu-time cost.

\begin{table}[!htbp]
	\centering
	\caption{Detailed information for results of the ABC flow with $\epsilon=0.01$ obtained by setting $\gamma_1=\gamma_2=1.0-1.0\times 10^{-8}$ and  $\gamma_3=0.9,0.99, 0.999, 0.9999$, respectively.}
	\label{Tab-ABC-0.01-gamma3}
	\begin{tabular}{| c| c|c | c | c |c|}
		\hlinew{1.2pt}
		&POD-updates&   DOFs   &   Average Error&  Times (s) \\
		\hlinew{1.2pt}
	$\gamma_3=0.9$   & 22 & 78  & 0.078965 & 19831.24 \\
$\gamma_3=0.99$    & 21 &  239& 0.043133 & 37962.03 \\
$\gamma_3=0.999$   & 19  & 393 &  0.015433 & 49428.57 \\
$\gamma_3=0.9999$ & 7 & 318 & 0.043027  &  18551.14  \\
		\hlinew{1.5pt}
	\end{tabular}
\end{table}

	\begin{table}[!htbp]
	\centering
	\caption{Detailed information for results of the ABC flow with $\epsilon=0.01$ obtained by setting $\gamma_3=1.0-1.0\times 10^{-8}$ and   $\gamma_1=\gamma_2=0.9,0.99, 0.999, 0.9999$, respectively.}
	\label{Tab-ABC-0.01-gamma1-2}
	\begin{tabular}{| c| c|c | c | c |c|}
		\hlinew{1.2pt}
		&POD-updates&   DOFs   &   Average Error&  Times (s) \\
		\hlinew{1.2pt}
			$\gamma_1=\gamma_2=0.9$   & 14 & 134 &0.037092  & 13922.37\\
$\gamma_1=\gamma_2=0.99$    & 7 & 146 & 0.011630 & 10635.61\\
$\gamma_1=\gamma_2=0.999$   & 6 & 192 &  0.015851 &  10279.38 \\
$\gamma_1=\gamma_2=0.9999$ & 6 & 237 &  0.020604 &   12319.63 \\
		\hlinew{1.5pt}
	\end{tabular}
\end{table}

\begin{table}[!htbp]
	\centering
	\caption{Detailed information for results of the ABC flow with $\epsilon=0.01$ obtained by setting  $\gamma_1=\gamma_2=\gamma_3=\gamma=1.0-1.0\times 10^{-2},1.0-1.0\times 10^{-4},1.0-1.0\times 10^{-6},1.0-1.0\times 10^{-8}$ respectively.}
	\label{Tab-ABC-0.01-gamma}
	\begin{tabular}{| c| c |c| c | c |c|}
		\hlinew{1.2pt}
		&  POD-updates	&POD-DOFs&  Average Error& Times (s) \\
		\hlinew{1.2pt}
	 			$\gamma=1-10^{-2}$   & 24  & 78 & 0.022985 & 19687.06 \\
	$\gamma=1-10^{-4}$    & 7 & 262   &0.014736    & 14718.30 \\
	$\gamma=1-10^{-6}$  & 6 & 275 & 0.031665 & 14256.98 \\
	$\gamma=1-10^{-8}$ & 5 & 245 & 0.054979  & 11788.47  \\
		\hlinew{1.5pt}
	\end{tabular}
\end{table}

\section{Numerical experiments for different coarse grid}

Here, we use some numerical experiments to show how the degree of freedom for the coarse grid affects the accuracy and the cpu time cost when the fine grid is fixed. We hope it can provide some useful information about how to choose the coarse grid. 

We  take the Kolmogorov flow with $\epsilon = 0.01$ as an example to show the behavior of our two-grid adaptive POD algorithm with different coarse mesh. 
In our experiments, we fix the fine mesh as what is shown in the manuscript, the other parameters are also set as the same as those used for obtaining results shown in Table 1. That is, we set $\gamma_1=\gamma_2 = 0.999$, $\gamma_3=1.0-1.0\times 10^{-8}$, and $\eta_0 = 0.005$. The detailed  numerical results are listed in Table \ref{table-A1}. 

In Table \ref{table-A1}, `N-Refine' means the number of refinement  used to obtain the mesh from the initial mesh, `Time step' means the time step used for the coarse mesh, `DOFs-CoarseGrid' means the degree of freedom corresponding to the coarse grid,  `DOFs-POD' means number of POD modes corresponding to the coarse grid, 
`POD-updates' means the number of updation of the POD modes, `Average Error' is computed by averaging the errors of numerical solution for each time instance, `Time-CoarseGrid' is the wall time for the simulation on the coarse grid, `Total-Time' is the wall time for all the simulation.

\begin{table}[!htbp]
	\centering
	\caption{The results of Kolmogorov flow with $\epsilon=0.01$ obtained by TG-APOD with different coarse mesh, respectively.}
	\label{table-A1}
	\begin{tabular}{|c|c|c|c|c|c|} \hlinew{1.5pt}	
		N-Refine & 15 &16 &17 &18&19 \\ 
		Time step &  0.1 & 0.09 & 0.05& 0.03& 0.01 \\ 
		DOFs-CoarseGrid &32768  & 65536 & 163840 & 262144 & 524288 \\ 
		Time-CoarseGrid & 17.47 & 34.03 & 133.46 & 882.47 & 1087.41 \\ \hlinew{1.5pt}
		POD-updates & 7 & 7 & 8 & 8 & 8 \\
		DOFs-POD & 155 &  138 & 175 & 174 & 175  \\
		Average Error & 0.011799 &  0.010969  & 0.006300 & 0.004592 & 0.003830 \\
		Total-Time(s) & 7350.98 & 6299.70 & 9787.82 & 9312.94 & 9526.88 \\ \hlinew{1.5pt}
	\end{tabular}
\end{table}

\begin{figure}[!htbp]	
	\centering
	\begin{minipage}[t]{0.8\textwidth} 
		\includegraphics[width =\textwidth]{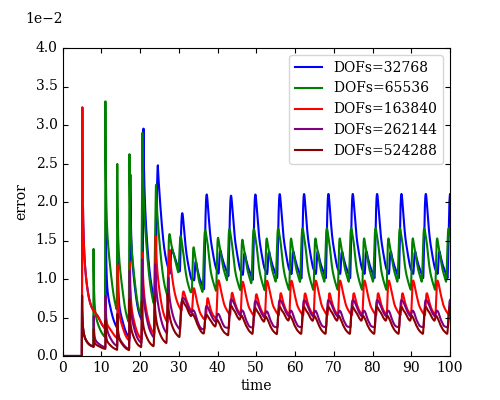}
	\end{minipage}
	\caption{The change of error for the solution of Kolmogorov flow with $\epsilon=0.01$ obtained by TG-APOD with different coarse mesh.}
	\label{coarse-mesh-compare}
\end{figure}

To see more clearly, we compare the error of numerical solutions obtained by TG-APOD with different coarse mesh in Fig \ref{coarse-mesh-compare}.

From Table \ref{table-A1} and Fig. \ref{coarse-mesh-compare}, we can see that as  the degree of freedom for the coarse mesh increases, the average error will decrease, while the cpu time cost may increase. Therefore, there is a balance for the accuracy and the cpu time cost. For our case, taking into account  the accuracy and the cpu time cost, we choose the coarse grid obtained by refine $16$ times from the initial mesh. 

If we take a detailed look at Table \ref{table-A1}, we will find that the cpu time cost is not only affected by the degree of freedom for the coarse grid. In fact, if the degree of freedom for the coarse mesh are not so large,  the cpu time will be mainly determined by the number of POD modes updates (the process of basis update requires discretize the original equation with finite element basis for a period of time, which is expensive), which is mainly determined by the $\eta_0$. However, as  the degree of freedom for the coarse mesh increases,  the cpu time cost by calculating the error indicator on the coarse mesh also needs to be considered. Anyway,  it is not easy to determine how large the coarse mesh is for the maximum efficiency.

\end{appendix}
  

  
\bibliographystyle{spmpsci}      
 
\bibliography{ref}
 
\end{document}